\documentclass{amsart}
\usepackage{graphicx}
\usepackage{amssymb}
\usepackage{epstopdf}
\usepackage{nicefrac}
\usepackage{enumerate}
\usepackage{hyperref}
\usepackage{float}
\usepackage{color}
\usepackage{pst-all}
\usepackage{tabularx}

\newtheorem{thm}{THEOREM}[section]
\newtheorem{conj}[thm]{CONJECTURE}
\newtheorem{cor}[thm]{COROLLARY}

\newtheorem{defn}[thm]{DEFINITION}

\newtheorem{lemma}[thm]{LEMMA}
\newtheorem{prob}[thm]{PROBLEM}
\newtheorem{prop}[thm]{PROPOSITION}
\newtheorem{quest}[thm]{QUESTION}
\newtheorem{remark}[thm]{REMARK}

\newcommand{\ds}{\displaystyle}

\newcommand{\F}{{\mathcal F}}
 
\newcommand{\FfM}{{\mathcal F}_{\fM}}

 \newcommand{\bH}{{\bf H}}
 \newcommand{\bK}{{\bf K}}
 
 \newcommand{\cA}{{\mathcal A}}
 \newcommand{\cB}{{\mathcal B}}
 \newcommand{\cC}{{\mathcal C}}
 \newcommand{\cD}{{\mathcal D}}
\newcommand{\cG}{{\mathcal G}}
\newcommand{\cH}{{\mathcal H}}
\newcommand{\cK}{{\mathcal K}}

\newcommand{\cO}{{\mathcal O}}
\newcommand{\cP}{{\mathcal P}}

\newcommand{\cS}{{\mathcal S}}
\newcommand{\cW}{{\mathcal W}}

\newcommand{\e}{{\varepsilon}} 
\newcommand{\G}{\Gamma}
\newcommand{\mR}{{\mathbb R}}
\newcommand{\mS}{{\mathbb S}}
\newcommand{\mT}{{\mathbb T}}

\newcommand{\oU}{{\overline U}}

\newcommand{\wtpi}{\widetilde{\pi}}

\newcommand{\mZ}{{\mathbb Z}}

\newcommand{\whG}{{\widehat{G}}}
\newcommand{\whH}{{\widehat{H}}}

\newcommand{\whU}{{\widehat U}}

\newcommand{\whW}{{\widehat W}}

\newcommand{\whPhi}{{\widehat{\Phi}}}
\newcommand{\fM}{{\mathfrak{M}}}

\newcommand{\fC}{{\mathfrak{C}}}

\newcommand{\fT}{{\mathfrak{T}}}
\newcommand{\fX}{{\mathfrak{X}}}

\newcommand{\fU}{{\mathfrak{U}}}

\newcommand{\wtM}{\widetilde{M}}

 \newcommand{\ovq}{\overline{q}}

\newcommand{\vp}{{\varphi}}
\newcommand{\whvarp}{{\widehat \varphi}}
\newcommand{\cV}{{\mathcal V}}

\newcommand{\cU}{{\mathcal U}}
\newcommand{\cT}{{\mathcal T}}

\newcommand{\cGF}{\cG_{\F}} 
\newcommand{\cI}{{\mathcal I}}

 \newcommand{\whL}{\widehat{L}}
\newcommand{\whe}{\widehat{e}}
 
 \newcommand{\whfM}{\widehat{\fM}}
 \newcommand{\whcS}{\widehat{\cS}}

 \newcommand{\whiota}{\widehat{\iota}}
 \newcommand{\whphi}{\widehat{\phi}}
 \newcommand{\whrho}{\widehat{\rho}}

\newcommand{\whx}{{\widehat{x}}}
\newcommand{\whd}{{\widehat{d}}}
\newcommand{\whg}{{\widehat{g}}}
\newcommand{\whh}{{\widehat{h}}}
\newcommand{\whk}{{\widehat{k}}}

\newcommand{\whq}{{\widehat{q}}}

\newcommand{\wtSigma}{{\widetilde{\Sigma}}}
\newcommand{\wtF}{{\widetilde{\F}}}
 
 \newcommand{\psg}{{\rm pseudo}{\star}{\rm group}}

\newcommand{\dX}{d_{\fX}} 
\newcommand{\dF}{d_{\F}} 
\newcommand{\dM}{d_{\fM}}

\newcommand{\fD}{{\mathfrak{D}}}
\newcommand{\fR}{{\mathfrak{R}}}

\newcommand{\ocH}{\overline{\mathcal H}}
\newcommand{\orho}{\overline{\rho}}

\input xy
\xyoption {all}

 \newcommand\tail{\mathrel{\stackrel{\makebox[0pt]{\mbox{\normalfont\tiny t}}}{\sim}}}

\begin{document}

\title{Wild solenoids}
 
\author{Steven Hurder}
 \author{Olga Lukina}
 \email{hurder@uic.edu, lukina@uic.edu}
\address{SH \& OL: Department of Mathematics, University of Illinois at Chicago, 322 SEO (m/c 249), 851 S. Morgan Street, Chicago, IL 60607-7045}

\thanks{Version date:  February 9, 2017}



\dedicatory{Dedicated to the memory of James T. Rogers, Jr.}

 \begin{abstract}
A weak solenoid is a foliated space defined as the inverse limit of finite coverings of a closed compact manifold $M$. The monodromy of a  weak solenoid defines an equicontinuous minimal action on a Cantor space $X$ by the fundamental group $G$ of $M$. The discriminant group of this action is an obstruction to this action being homogeneous. The discriminant vanishes if the group $G$ is abelian, but there are examples of actions of nilpotent groups  for which the discriminant is non-trivial.  The action is said to be stable if the discriminant group remains unchanged for  the induced action on sufficiently small  clopen neighborhoods in $X$. If the discriminant group never stabilizes as the diameter of the clopen set $U$ tends to zero, then we say that the action  is unstable, and the weak solenoid which defines it is said to be wild. In this work, we show two main results in the course of our study of the properties of the discriminant group for Cantor actions. First, the tail equivalence class of the sequence of discriminant groups obtained for the restricted action on a neighborhood basis system of a point in $X$ defines an invariant of the return equivalence class of the action, called the asymptotic discriminant, which is consequently an invariant of the homeomorphism class of the weak solenoid. Second, we construct uncountable collections of wild solenoids with pairwise distinct asymptotic discriminant invariants for a fixed base manifold $M$, and hence fixed finitely-presented group $G$, which are thus pairwise non-homeomorphic. The study in this work is the continuation of the seminal works on homeomorphisms of weak solenoids by   Rogers and Tollefson in 1971, and is dedicated to the memory of Jim Rogers. 
\end{abstract}

\maketitle

 \vfill
 \eject

\section{Introduction}\label{sec-intro}

A weak solenoid $\cS_\cP$, as introduced by McCord in \cite{McCord1965}, is the inverse limit space of a sequence of proper covering maps 
  $$\cP = \{p_{\ell +1} \colon M_{\ell +1} \to M_{\ell} \mid \ell \geq 0\},$$ 
  where $M_\ell$ is a compact connected manifold without boundary, and $p_{\ell+1}$ is a finite-to-one covering map. The collection of maps $\cP$ is called a \emph{presentation} for $\cS_\cP$.  A solenoid $\cS_\cP$ is \emph{regular} if the presentation $\cP$ can be chosen so that for each $\ell \geq 0$ the composition $p_\ell \circ \cdots p_1:M_\ell \to M_0$ is a regular covering map.  A weak solenoid which is not regular is said to be \emph{irregular}, and the purpose of this work is to study the classification problem for irregular solenoids.

Let $\fM$ be an \emph{$n$-dimensional foliated space} $\fM$ with foliation $\FfM$, as   defined   in Section~\ref{sec-concepts}. A foliated space is a natural generalization of the notion of a smooth foliation $\F$ of a compact connected manifold $M$. 
A foliated space $\fM$ for which the local transversals to the foliation $\FfM$ are totally disconnected is called a \emph{matchbox manifold}, and   the leaves of $\FfM$ are the path connected components of $\fM$.  Matchbox manifolds were introduced in the $1$-dimensional case in the works \cite{AM1988,AO1991,AO1995}, have been studied in the higher dimensional cases in the works     \cite{ClarkHurder2011,ClarkHurder2013,CHL2013a,CHL2014,CHL2013c,DHL2016a,DHL2016b,DHL2016c}. 

A space is said to be  \emph{foliated homogeneous} if given any pair of points $x,y \in \fM$, then there exists a leaf-preserving homeomorphism $h \colon \fM \to \fM$ such that $h(x) = y$. 
A simple observation about matchbox manifolds is that a homeomorphism   $\varphi \colon \fM \to \fM'$    preserves the path-connected components, hence  a homeomorphism of a matchbox manifolds preserves their foliations. It follows that if $\fM$ is   homogeneous, then it is   foliated homogeneous.

 McCord observed in \cite{McCord1965} that  a weak solenoid is  an equicontinuous matchbox manifold with minimal leaves and equicontinuous holonomy pseudogroup.   The following result was obtained in \cite{Thomas1973} for the $1$-dimensional case, and in  \cite{ClarkHurder2013} for matchbox manifolds with leaves of dimension $2$ and higher.

\begin{thm}\cite{ClarkHurder2013}\label{thm-hmm}
Let $\fM$ be a homogeneous matchbox manifold. Then:
\begin{enumerate}
\item The holonomy pseudogroup of $\FfM$ is equicontinuous;
\item $\FfM$ is minimal, that is, every leaf of $\FfM$ is dense in $\fM$;
\item $\fM$ is homeomorphic to a regular   solenoid.
\end{enumerate}
\end{thm}

\begin{thm}\label{thm-mmweaksolenoid}\cite{ClarkHurder2013}
Let $\fM$ be a matchbox manifold. Then the holonomy pseudogroup of $\FfM$ is equicontinuous if and only if there exists a homeomorphism $\phi \colon \fM \to \cS_\cP$ for a presentation $\cP$.   The solenoid $\cS_\cP$ is regular if and only if $\fM$ is homogeneous. 
\end{thm}

The construction of the homeomorphism between $\fM$ and a weak solenoid $\cS_\cP$ depends on the choice of a clopen transversal to $\FfM$. Given  two such choices of clopen transversals to $\FfM$, the resulting weak solenoids are related by return equivalence \cite{CHL2013c}, which is discussed in Section~\ref{subsec-return}. A key result, as stated in Theorem~\ref{thm-topinv}, is that homeomorphic matchbox manifolds have return equivalent holonomy pseudogroups.

Thus, a basic problem for the study of the classification of weak solenoids up to homeomorphism, is to develop invariants for the return equivalence class of a minimal pseudogroup action on a Cantor space.
In Section~\ref{sec-adi}, we define the \emph{asymptotic discriminant}   for weak solenoids, and show that it is an   invariant of the homeomorphism type for equicontinuous matchbox manifolds.

The results of  this work are also closely related to another concept which was introduced in the study of smooth Riemannian foliations. 

\begin{defn}\label{molinospace}
A \emph{Molino space} for $(\fM, \FfM)$ is a (connected) foliated homogeneous space   $(\whfM, \widehat{\FfM})$ of the same leaf dimension as $\fM$, and a foliated principal $\bH$-fibration 
 \begin{equation}\label{eq-molinoseqM}
\bH  \longrightarrow \whfM \stackrel{\whq}{\longrightarrow} \fM \ , 
\end{equation}
 where $\bH$ is a compact topological group,   $\whfM$ is foliated homogeneous, and for each leaf $\whL \subset \whfM$ there is a leaf $L \subset \fM$ such that the restriction $\whq | \whL \to L$ is a covering map.  
\end{defn}
We call  \eqref{eq-molinoseqM} a \emph{Molino sequence} for the foliation.

  Molino showed in \cite{Molino1982,Molino1988} (see also \cite{Haefliger1989,MoerdijkMrcunbook2003})  that a Riemannian foliation $\F$ of a smooth compact manifold $M$ always admits a Molino sequence. In this case, the fiber group $\bH$ is a connected compact Lie group whose Lie algebra is isomorphic to  the Lie algebra of transverse projectable vector fields for $\F$. Thus, $\bH$ is locally defined, up to local  isomorphism, by the transverse geometry to $\F$. 
 
 In a series of papers, \cite{ALC2009,ALC2010,ALM2016,ALB2016},  {\'A}lvarez L{\'o}pez and his coauthors   formulated a \emph{topological Molino theory} for equicontinuous foliated spaces, which is a partial generalization of  the Molino theory for smooth Riemannian foliations.   They formulated the notion of  \emph{strongly quasi-analytic regularity} for the holonomy pseudogroup of a  foliated space, or the \emph{SQA condition}, which is   discussed in  Section~\ref{subsec-defSQA}.  The topological Molino theory developed in \cite{ALM2016} applies to foliated spaces which satisfy the SQA condition.

A Molino space for a matchbox manifold $\fM$ is a homogeneous  matchbox manifold $\whfM$, which is an equicontinuous foliated space by Theorem~\ref{thm-hmm}. The quotient foliated space of an equicontinuous foliated space is equicontinuous, so by Theorems~\ref{thm-hmm} and \ref{thm-mmweaksolenoid} we have:
\begin{cor}\label{cor-weaksols}
Let $\fM$ be a matchbox manifold which admits a Molino space $\whfM$ with Molino sequence \eqref{eq-molinoseqM}. Then $\fM$ is homeomorphic to a weak solenoid.
\end{cor}
Thus, the study of Molino spaces and their properties for matchbox manifolds reduces to the study of weak solenoids. A converse to Corollary \ref{cor-weaksols} was obtained by the authors joint with Dyer in \cite{DHL2016c}.

  \begin{thm}\cite[Theorem~1.2]{DHL2016c}\label{thm-molinospaces}
Let     $\fM$ be an equicontinuous matchbox manifold, and let $\cP$ be a presentation with homeomorphism $\fM \to \cS_\cP$ to a weak solenoid $\cS_\cP$.  
Then there exists a   regular  solenoid $\whcS_{\cP}$, a  compact totally disconnected  group $\cD_{\cP}$, and a  Molino sequence
\begin{equation}\label{eq-molinoseq}
\cD_{\cP} \longrightarrow \whcS_{\cP} \stackrel{\whq}{\longrightarrow} \cS_{\cP} \ ,
\end{equation}
where the    spaces $\whcS_{\cP}$ and $\cD_{\cP}$ depend on the choice of the presentation $\cP$. 
  \end{thm}

One motivation for  the work in this paper is to investigate how the the space $\whcS_{\cP}$ depends on the choice of the presentation $\cP$. In particular, \cite{DHL2016c} describes examples where the cardinality of the fibre group $\cD_\cP$ varies with the choice of the presentation $\cP$. 
So it is important to understand when the resulting Molino sequence \eqref{eq-molinoseq} is well-defined, up to topological conjugacy of fibrations. This is equivalent to finding conditions such that   the   homeomorphism type of the fiber group $\cD_{\cP}$ is well-defined,   that is, it is independent of the choice of the presentation for $\cS_{\cP}$ and hence is an invariant of the homeomorphism type of $\fM$.

\begin{defn}\label{def-stable}
An equicontinuous matchbox manifold $\fM$ is said to be \emph{stable} if its associated Molino sequence \eqref{eq-molinoseq} is well-defined, up to topological conjugacy of fibrations, and otherwise is said to be \emph{wild}.
\end{defn}

The construction of a Molino sequence for a weak solenoid $\cS_{\cP}$  was given in \cite{DHL2016c}, using the methods of group chains, which we recall in Section~\ref{sec-stablechains}. Moreover, a criteria in terms of group chains was given in that work for when $\cS_{\cP}$ is stable. The work \cite{DHL2016c} also proved stability for a few specific classes of weak solenoids. Remarkably, it turned out that the holonomy pseudogroups of solenoids in these classes satisfies the \emph{strong quasi-analytic property} (SQA condition), introduced by \'Alvarez L\'opez and Candel \cite{ALC2009}. This property is an analog for foliated spaces of the \emph{quasi-analytic property} for real-analytic foliations and it is described in detail in Section~\ref{subsec-defSQA}. 

  The first main theorem of this paper shows that there is a direct link between the SQA property of \'Alvarez L\'opez and Candel \cite{ALC2009} and the stability of Molino sequence, as defined in \cite{DHL2016c}. 

\begin{thm}\label{thm-main1}
An equicontinuous matchbox manifold $\fM$  is stable if and only if its holonomy pseudogroup satisfies the    LCSQA  (Locally  Completely Strongly Quasi-Analytic) property introduced in Definition~\ref{def-LCSQA}.
\end{thm}
Proposition~\ref{prop-sqaequiv} gives the proof of this result, and Proposition \ref{prop-retequivstable} shows the intuitively clear result that the LCSQA property is an invariant of return equivalence for pseudogroup actions. 

 The \emph{Schori solenoid} was constructed in \cite{Schori1966} as one of the first examples of non-homogeneous solenoids. The Schori solenoid is the inverse limit of $3$-to-$1$ coverings of the genus $2$ surface. It was proved in \cite{DHL2016c} that the holonomy pseudogroup of this solenoid is not LCSQA, and thus by  Theorem~\ref{thm-main1}, we obtain that it provides the first example of a wild solenoid.
\begin{cor}
The Schori solenoid is wild.
\end{cor}

 Another result of this paper is the introduction of an invariant, which allows us to distinguish between various classes wild solenoids. This invariant, called the \emph{the asymptotic discriminant}, is defined in Section~\ref{sec-adi}, where it is also shown to be an invariant under return equivalence. The results of this section combine  show the following:

  \begin{thm}\label{thm-AD1}
 The  asymptotic discriminant  of an equicontinuous matchbox manifold $\fM$ is well-defined, and defines an invariant of  the homeomorphism class of $\fM$.
  \end{thm}

Thus, the natural question is to ask whether this invariant is effective, and how does one construct examples which are distinguished by the asymptotic discriminant? 
  The answer to this question occupies the remainder of the paper, beginning in Sections~\ref{sec-solenoids} and \ref{sec-stablechains}. 
 Section~\ref{sec-construction} then gives a method for constructing weak solenoids which are wild and have prescribed asymptotic discriminant invariants. We apply this method in Section~\ref{sec-wildex} to show the following:

\begin{thm}\label{thm-main2}
For $n \geq 3$, let $G \subset {\bf SL}(n,\mZ)$ be a torsion-free subgroup of finite index. Then there exists uncountably many distinct homeomorphism types of  weak solenoids which are wild, all with the same base manifold $M_0$ having fundamental group $G$. 
\end{thm}

It is interesting to compare this result with some previous results about the classification problem for solenoids. 
The classification of $1$-dimensional solenoids up to homeomorphism was completed by  Bing \cite{Bing1960} and McCord \cite{McCord1965}, who showed that the homeomorphism type of the inverse limit $\cS_{\cP}$ of a $1$-dimensional presentation $\ds \cP = \{ p_{\ell+1} \colon \mS^1 \to \mS^1 \mid \ell \geq 0\}$ is determined by the   collection of the covering degrees
$\cD_k = \{ d_{\ell} = \deg\{ p_{\ell}\} \geq 2 \mid \ell \geq k\}$ for $k$ large. This classification was reformulated in terms of return equivalence for the solenoid $\cS_{\cP}$ by Aarts and Fokkink    \cite{AartsFokkink1991}; see also  \cite[Section~8.1]{CHL2013c}.

For a weak solenoid $\cS_{\cP}$  with leaf dimension $n \geq 2$,  the authors joint with Alex Clark obtained the following   result:

\begin{thm}\cite[Theorem~1.4]{CHL2013c}\label{thm-topinv3}
Let $\cS_{\cP}$ and $\cS'_{\cP'}$ be  weak   solenoids which are both $\mT^n$-like.  Then $\cS_{\cP}$ and $\cS'_{\cP'}$  are homeomorphic if and only if  their global monodromy actions $(\fX, G, \Phi)$ and $(\fX', G', \Phi')$ are return equivalent, where we have $G \cong G' \cong \mZ^n$.
\end{thm}

In general, for a weak solenoid $\cS_{\cP}$ with leaf dimension $n \geq 2$,
the fundamental group $G_0 = \pi_1(M_0, x_0)$ of the base manifold $M_0$ of the presentation $\cP$ need not be abelian, and the classification problem is far from being well-understood.    By Theorem~\ref{thm-AD} and Proposition~\ref{prop-ADequivalence}, 
the asymptotic discriminant   associated to a equicontinuous matchbox manifold  provides a new invariant which can be used to show that a pair, or even an infinite family of solenoids  are non-homeomorphic.

We conclude with an observation about the embedding property for weak solenoids. 
It follows from the definitions that if a weak solenoid $\cS_{\cP}$ is homeomorphic to a minimal set for a $C^r$-foliation $\F$ of a manifold, in the sense of the works \cite{ClarkHurder2011,Hurder2016}, then in the real analytic case where  $r=\omega$, the holonomy  pseudogroup of $\cS_{\cP}$ satisfies the LCSQA condition.   We thus  obtain a  non-embedding result:
\begin{cor}\label{cor-realization}
The wild solenoids $\cS_{\cP}$ constructed in the proofs of Theorem~\ref{thm-main2} are not homeomorphic to the minimal sets of a real-analytic foliation on a manifold. 
\end{cor}

A solution to the following problem would represent an extension of the ideas in \cite{Ghys1993}:
\begin{conj}\label{conj-realization}
A wild solenoid $\cS_{\cP}$ is not homeomorphic to the minimal set of any $C^2$-foliation of a finite dimensional manifold.
\end{conj}

\section{Foliated spaces and return equivalence} \label{sec-concepts}

In this section, we briefly recall   the definitions of foliated spaces and matchbox manifolds, and some of their  attributes.  In particular, we recall  in Theorem~\ref{thm-topinv}    that a homeomorphism between minimal matchbox manifolds induces return equivalence of their holonomy pseudogroups. More details and properties of foliated spaces  can be found in   the works \cite{CandelConlon2000,ClarkHurder2013,CHL2013a,CHL2014,MS2006}.  

\subsection{Foliated spaces}\label{subsec-fs}
Recall that a \emph{continuum} is a compact connected metrizable space.
\begin{defn} \label{def-fs}
A \emph{foliated space of dimension $n$} is a   continuum $\fM$, such that  there exists a compact separable metric space $\fX$, and
for each $x \in \fM$ there is a compact subset $\fT_x \subset \fX$, an open subset $U_x \subset \fM$, and a homeomorphism defined on the closure
$\vp_x \colon \oU_x \to [-1,1]^n \times \fT_x$ such that $\vp_x(x) = (0, w_x)$ where $w_x \in int(\fT_x)$. Moreover, it is assumed that each $\vp_x$  admits an extension to a foliated homeomorphism
$\whvarp_x \colon \whU_x \to (-2,2)^n \times \fT_x$ where $\whU_x$ is an open subset such that $\oU_x \subset \whU_x$.
\end{defn}
The subspace $\fT_x$ of $\fX$ is  called the \emph{local transverse model} at $x$. In the case where  $\fM$ is a connected topological manifold, we say that $\fM$ is a (topological) foliated manifold. Then  the maps $\vp_x$ are coordinate charts for $\fM$,  and the sets $\fT_x $ are connected open subsets of some Euclidean space $\mR^q$. At the other extreme,  a \emph{matchbox manifold} is a    foliated space $\fM$ such that for each $x \in \fM$, the transverse model space $\fT_x \subset \fX$ is totally disconnected.

Let $\pi_x \colon \oU_x \to \fT_x$ denote the composition of $\vp_x$ with projection onto the second factor.

For $w \in \fT_x$ the set $\cP_x(w) = \pi_x^{-1}(w) \subset \oU_x$ is called a \emph{plaque} for the coordinate chart $\vp_x$. We adopt the notation, for $z \in \oU_x$, that $\cP_x(z) = \cP_x(\pi_x(z))$, so that $z \in \cP_x(z)$. Note that each plaque $\cP_x(w)$ is given the topology so that the restriction $\vp_x \colon \cP_x(w) \to [-1,1]^n \times \{w\}$ is a homeomorphism. Then $int (\cP_x(w)) = \vp_x^{-1}((-1,1)^n \times \{w\})$.
Let $U_x = int (\oU_x) = \vp_x^{-1}((-1,1)^n \times int(\fT_x))$.
Note that if $z \in U_x \cap U_y$, then $int(\cP_x(z)) \cap int( \cP_y(z))$ is an open subset of both
$\cP_x(z) $ and $\cP_y(z)$.
The collection
$$\cV = \{ \vp_x^{-1}(V \times \{w\}) \mid x \in \fM ~, ~ w \in \fT_x ~, ~ V \subset (-1,1)^n ~ {\rm open}\}$$
forms the basis for the \emph{fine topology} of $\fM$. The connected components of the fine topology are called leaves, and define the foliation $\FfM$ of $\fM$.
For a matchbox manifold,  the leaves of  $\FfM$  are the  maximal path-connected components of $\fM$.
Let $L_x \subset \fM$ denote the leaf of $\FfM$ containing $x \in \fM$.

A map $h \colon \fM \to \fM'$ between foliated spaces is said to be a \emph{foliated map} if the image of each leaf of $\FfM$ is contained in a leaf of $\F_{\fM'}$.  Each leaf of $\FfM$ is path connected, so its image under $h$ is path connected. If  $\fM'$ is a matchbox manifold, then the image of a leaf of $\FfM$ \emph{must} be contained in a leaf of $\F_{\fM'}$. Thus we have:
\begin{lemma} \label{lem-foliated1}
Let $\fM$ and $\fM'$ be matchbox manifolds, and $h \colon \fM \to \fM'$ a continuous map. Then $h$ maps the leaves of $\FfM$ into leaves of $\F_{\fM'}$. In  particular, any homeomorphism $h \colon \fM \to \fM'$ of   matchbox manifolds is a foliated map.  
\end{lemma}

A   foliated space $\fM$ is said to be \emph{smooth}, if the local charts $\vp_x \colon \oU_x \to [-1,1]^n \times \fT_x$ can be chosen so that for each leaf $L \subset \fM$ of $\FfM$,  the restrictions of the charts to $L$   defines a smooth manifold structure on $L$.
That is, for all $x,y \in \fM$ with $z \in U_x \cap U_y$, there exists an open set $z \in V_z \subset U_x \cap U_y$ such that $\cP_x(z) \cap V_z$ and $\cP_y(z) \cap V_z$ are connected open sets, and the composition
$$\psi_{x,y;z} \equiv \vp_y \circ \vp_x ^{-1}\colon \vp_x(\cP_x (z) \cap V_z) \to \vp_y(\cP_y (z) \cap V_z)$$
is a smooth map, where $\vp_x(\cP_x (z) \cap V_z) \subset \mR^n \times \{w\} \cong \mR^n$ and $\vp_y(\cP_y (z) \cap V_z) \subset \mR^n \times \{w'\} \cong \mR^n$. The leafwise transition maps $\psi_{x,y;z}$ are assumed to depend continuously on $z$ in the $C^{\infty}$-topology.

We recall a standard result, whose proof for foliated spaces can be found in \cite[Theorem~11.4.3]{CandelConlon2000}.
\begin{thm}\label{thm-riemannian}
Let $\fM$ be a smooth foliated space. Then there exists a leafwise Riemannian metric for $\FfM$, such that for each $x \in \fM$, the leaf $L_x$ inherits the structure of a complete Riemannian manifold with bounded geometry, and the Riemannian geometry depends continuously on $x$ .  
\end{thm}
All matchbox manifolds  are assumed to be smooth with a fixed choice of    metric $\dM$ on $\fM$, and a choice of a compatible    leafwise Riemannian metric $\dF$.

\subsection{Holonomy pseudogroup} \label{subsec-holonomy}

We next describe the construction of the holonomy pseudogroup associated to a matchbox manifold. We only give details sufficient for defining the notion of return equivalence in Section~\ref{subsec-return} below, and refer to the text   \cite{CandelConlon2000} and the works \cite{ClarkHurder2013,CHL2014,CHL2013c} for   details. 

The first step is to choose a     \emph{regular foliated covering}  of $\fM$,   $\ds \{\vp_i \colon  U_i \to (-1,1)^n \times \fT_i \mid 1 \leq i \leq \nu\}$ where 
$\cU = \{U_{1}, \ldots , U_{\nu}\}$ is an  open covering of $\fM$, and each chart $\vp_i \colon  U_i \to (-1,1)^n \times \fT_i$ satisfies the conditions on charts in  Definition~\ref{def-fs}, and also 
  satisfies the regularity conditions detailed in Proposition~2.6 of  \cite{CHL2014}, which in particular includes the assumption that the plaques $\cP_x(z)$ associated to the foliation charts are convex in the leafwise metric.

We   also  assume that the transverse model spaces $\fT_i$ form a  disjoint  clopen covering  of   $\fX$, so that 
  $\ds \fX = \fT_1 \ \dot{\cup} \cdots \dot{\cup} \ \fT_{\nu}$.
For each $1 \leq i \leq \nu$, the set $\cT_i =  \vp_i^{-1}(0 , \fT_i)$ is a compact transversal to $\FfM$. Without loss of generality, we can assume   that the transversals 
$\ds \{ \cT_{1} , \ldots , \cT_{\nu} \}$ are pairwise disjoint in $\fM$, and let   $\cT = \cT_1 \cup \cdots \cup \cT_{\nu} \subset \fM$ denote their disjoint union. 
Define the maps $\ds \tau_i \colon \fT_i \to \oU_i$ given by  $\ds \tau_i(\xi) = \vp_i^{-1}(0 , \xi)$,  and let  $\tau \colon \fX \to \cT$ denote the union of the maps $\tau_i$. Each transversal $\cT_i$ has a metric $d_{\cT_i}$ obtained by the restriction of $\dM$, and we define 
   $$d_{\fT_i} (x,y) = \dM(\tau_i(x), \tau_i(y)) ~ {\rm if}~ x,y \in \fT_i, ~ {\rm for}~ 1 \leq i \leq \nu.$$
We may assume that $\dM$ is rescaled so that ${\rm diam}(\fM) <1$. The metric $\dX$ on $\fX$ is defined by setting 
$$\dX(x,y) =d_{\fT_i} (x,y) ~ {\rm if}~ x,y \in \fT_i \quad , \quad \dX(x,y) =  1 ~ {\rm otherwise} \ .$$

A \emph{leafwise path}  is a continuous map $\gamma \colon [0,1] \to \fM$ such that there is a leaf $L$ of $\FfM$ for which $\gamma(t) \in L$ for all $0 \leq t \leq 1$.
If $\fM$ is a matchbox manifold, and $\gamma \colon [0,1] \to \fM$ is continuous, then   $\gamma$ is a leafwise path by Lemma~\ref{lem-foliated1}. 
The holonomy pseudogroup for a matchbox manifold $(\fM, \FfM)$ is defined analogously to the construction of holonomy for foliated manifolds, and associates to a leafwise path with endpoints in the transversal space $\cT$ a local homeomorphism of $\fX$. In this way,  the construction  generalizes the induced dynamical systems  associated to a section of a flow. 

 For each  $1 \leq i \leq \nu$, let   $\pi_i  \colon \oU_i \to \fT_i$ be composition of the coordinate chart $\vp_i$ with the projection onto the second factor $\fT_i$.
A pair of indices $(i,j)$, $1 \leq i,j \leq \nu$, is said to be \emph{admissible} if the \emph{open} coordinate charts satisfy $U_i \cap U_j \ne \emptyset$.
For $(i,j)$ admissible, define clopen subsets $\fD_{i,j} = \pi_i(U_i \cap U_j) \subset \fT_i \subset \fX$.  The convexity of foliation charts imply that plaques are either disjoint, or have connected intersection. This implies that there is a well-defined homeomorphism $h_{j,i} \colon \fD_{i,j} \to \fD_{j,i}$ with domain $D(h_{j,i}) = \fD_{i,j}$ and range $R(h_{j,i}) = \fD_{j,i}$.

The maps $\cGF^{(1)} = \{h_{j,i} \mid (i,j) ~{\rm admissible}\}$ are the transverse change of coordinates defined by the foliation charts. By definition they satisfy $h_{i,i} = Id$, $h_{i,j}^{-1} = h_{j,i}$, and if $U_i \cap U_j\cap U_k \ne \emptyset$ then $h_{k,j} \circ h_{j,i} = h_{k,i}$ on their common domain of definition. The \emph{holonomy pseudogroup} $\cGF$ of $\FfM$ is the topological pseudogroup acting on $\fX$ generated by   the elements of $\cGF^{(1)}$ as follows. 

A sequence $\cI = (i_0, i_1, \ldots , i_{\alpha})$ is \emph{admissible}, if each pair $(i_{\ell -1}, i_{\ell})$ is admissible for $1 \leq \ell \leq \alpha$, and the composition
\begin{equation}\label{eq-defholo}
h_{\cI} = h_{i_{\alpha}, i_{\alpha-1}} \circ \cdots \circ h_{i_1, i_0}
\end{equation}
has non-empty domain.
The domain $\fD_{\cI}$ of $h_{\cI}$ is the \emph{maximal clopen subset} of $\fD_{i_0 , i_1} \subset \fT_{i_0}$ for which the compositions are defined.
Given any open subset $U \subset \fD_{\cI}$ define $h_{\cI} | U \in \cGF$ by restriction. 

The holonomy pseudogroup $\cGF$ is \emph{finitely-generated} if the generating set $\cGF^{(1)}$ is finite. In particular, if a foliated space $\fM$ is compact, then $\cGF$ is finitely-generated.
   
 \subsection{Pseudogroup dynamics}\label{subsec-psgdynamics}
 
Let $\Phi \colon \cG \times \fX \to \fX$ be a finitely-generated pseudogroup action, such as the pseudogroup $\cGF$ constructed above for a matchbox manifold. We assume there is a disjoint decomposition into clopen subsets of $\fX = \fT_1 \cup \cdots \cup \fT_{\nu}$, and a finite collection of maps $\cG^{(1)}  = \{h_{j,i} \mid (i,j) ~{\rm admissible}\}$ where  $\fD_{i,j} \subset \fT_i$ is a non-empty clopen subset,  $h_{i,j} \colon \fD_{i,j} \to \fD_{j,i}$ is a homeomorphism with $h_{j,i} = h_{i,j}^{-1}$, and $h_{i,i}$ is the identity map on $\fT_i$. Then $\cG$ is the pseudogroup  generated by  $\cG^{(1)}$. 
For  $g \in \cG$,  let $\fD(g) \subset \fX$ denote its domain, and $\fR(g) = g(\fD(g)) \subset \fX$ its range.

Define the maps $h_{\cI}$ as in \eqref{eq-defholo} for the admissible sequences $\cI = (i_0, i_1, \ldots , i_{\alpha})$.
Then  $\fD(h_{\cI}) \subset \fD_{i_0 , i_1}$ and $\fR(h_{\cI})   \subset \fD_{i_{\alpha}, i_{\alpha-1}}$.
 Also, introduce   the collection of maps:
     \begin{equation}\label{eq-restrictedgroupoid}
\cG^* = \left\{ h_{\cI} | U \mid \cI ~ {\rm admissible~ and}  ~ U \subset \fD(h_{\cI}) \right\} \subset \cG ~ .
\end{equation}
  The collection of maps $\cG^*$    is closed under the operations of compositions, taking inverses, and restrictions to open sets, and  is called a \emph{$\psg$} in the literature \cite{Hurder2014,Matsumoto2010}.

The $\cG$-orbit of a point  $w \in \fX$   is denoted by
\begin{equation}\label{eq-orbits}
\cO(w) = \{g ( w )\mid g \in \cG  ~ \textrm{such that}~ w \in \fD(g) \} = \{g ( w )\mid g \in \cG^*  ~\textrm{such that}~ w \in \fD(g) \} \subset \fX ~ .
\end{equation}
 
The action of $\cG$ is \emph{minimal} if for each $w \in \fX$ the $\cG$-orbit $\cO(w)$ is dense in $\fX$.

Let   $W \subset   \fX$ be an open subset, then    the \emph{restriction}   of $\cG^*$ to $W$ is defined by:
   \begin{equation}\label{eq-inducedpseudo}
\cG^*_W  =  \left\{ g \in \cG^* \mid \fD(g) \subset W   ~, ~ \fR(g) \subset  W   \right\} .
\end{equation}
The $\psg$ $\cG^*_W$  is   also referred to as the \emph{induced} $\psg$ on $W$.

 Introduce the filtrations of $\cG^*$   by word length.
For $\alpha \geq 1$, let  $\cG^{(\alpha)}$ be the collection of holonomy homeomorphisms $h_{\cI} | U \in \cG^*$ determined   by admissible paths $\cI = (i_0,\ldots,i_k)$ such that $k \leq \alpha$ and $U \subset \fD(h_{\cI})$ is open.
Then for  each $g \in \cG^*$ there is some $\alpha$ such that $g \in \cG^{(\alpha)}$. Let $\|g\|$ denote the least such $\alpha$, which is called the \emph{word length} of $g$.  

We   have the following finiteness result for minimal pseudogroup actions on compact spaces, whose counterpart for group actions is a well-known folklore result:
 \begin{lemma}\cite[Lemma~4.1]{CHL2014}\label{lem-finitegen}
Let $U \subset \fX$ be an open subset. Then there exists an integer  $\alpha_U$ such that $\fX$ is covered by the collection
$\{h_{\cI} (U)  \mid   h_{\cI} \in \cG^{(\alpha_U)} \}$, where  $h_{\cI} (U)  = h_{\cI} (U \cap \fD(h_{\cI}))$.
\end{lemma}

Next,      recall the     definition of an   equicontinuous  pseudogroup. Let $\dX$ denote the metric on $\fX$. 
\begin{defn} \label{def-equicontinuous}
The action of the    pseudogroup $\cG$ on $\fX$ is  \emph{equicontinuous} if for all $\epsilon > 0$, there exists $\delta > 0$ such that for all $g \in \cG^*$, if $w, w' \in \fD(g)$ and $\dX(w,w') < \delta$, then $\dX(g(w), g(w')) < \epsilon$.
That is, $\cG^*$ is equicontinuous as a family of local group actions.
\end{defn}

 Also, we define a notion of an analytic pseudogroup as follows:

\begin{defn} \label{def-analytic}
The action of the    pseudogroup $\cG$ on $\fX$ is  \emph{analytic} if   there exists an integer $q \geq 1$, an embedding $\fX \subset \mR^q$, and for all admissible pairs of indices $(i,j)$, open sets $U_{i,j} \subset \mR^q$ such that $\fD_{i,j}   \subset U_{i,j}$ and  map  $h_{j,i} \colon \fD_{i,j} \to \fD_{j,i}$  is the restriction of a real analytic map 
 $\whh_{j,i} \colon U_{i,j}  \to U_{j,i}$.
 \end{defn}

Finally, we recall a basic result from \cite{ClarkHurder2013}:
\begin{thm}\cite[Theorem~4.12]{ClarkHurder2013} \label{thm-minimal}
Let  $\cGF$ be the holonomy pseudogroup associated to a matchbox manifold $\fM$. If the  action of the    pseudogroup $\cGF$ on $\fX$ is  \emph{equicontinuous}, then the action is minimal. 
\end{thm}
The assumption that $\cGF$ is the holonomy pseudogroup derived from a connected foliated space $\fM$ is essential for the proof;  
the conclusion of Theorem~\ref{thm-minimal} is obviously false for an arbitrary pseudogroup.  

\subsection{Return equivalence}\label{subsec-return}

 The  concept of \emph{return equivalence} between  pseudogroups is the analog for matchbox manifolds of the   notion of \emph{Morita equivalence} for foliation groupoids associated to foliated manifolds, as discussed for example  in \cite{Haefliger1984}.

 For $i =1,2$, let    $\Phi^i \colon \cG_i \times \fX_i \to \fX_i$ be a finitely-generated pseudogroup action, as defined in Section~\ref{subsec-psgdynamics}.

 Given a clopen  set  $U_i \subset \fX_i$, for $i=1,2$,   let  $\cG^*_{U_i}$ denote the $\psg$ defined by the restriction of $\cG^*_i$ to $U_i$ as in \eqref{eq-inducedpseudo}.
We say that   $\cG^*_{U_1}$ and $\cG^*_{U_2}$
 are \emph{isomorphic} if there exists a homeomorphism $\phi \colon U_1 \to U_2$ such that the induced map $\Phi_{\phi} \colon \cG^*_{U_1} \to \cG^*_{U_2}$ is an isomorphism. That is,  for all $g \in \cG^*_{U_1}$ the map $\Phi_{\phi}(g) = \phi \circ g \circ \phi^{-1}$ defines an element of  $\cG^*_{U_2}$. Conversely, for   $h \in \cG^*_{U_2}$  the map $\Phi_{\phi}^{-1}(h) = \phi^{-1} \circ h \circ \phi$ defines an element of  $\cG^*_{U_1}$.

\begin{defn} \label{def-return}
Let    $\Phi^i \colon \cG_i \times \fX_i \to \fX_i$ be   minimal pseudogroup actions,  for $i=1,2$. Then $\Phi^1$ and $\Phi^2$ are \emph{return equivalent} if there are non-empty open sets $U_i \subset \fX_i$ and a homeomorphism $\phi \colon U_1 \to U_2$ such that the induced map $\Phi_{\phi} \colon \cG^*_{U_1} \to \cG^*_{U_2}$ is an isomorphism. 
 \end{defn}

We first note the following result:

\begin{prop}\cite[Proposition~4.6]{CHL2013c} \label{prop-return}
\emph{Return equivalence} is an equivalence relation on the class of finitely-generated minimal pseudogroup actions.
\end{prop}

The notion of return equivalence is applied to matchbox manifolds as follows:  

 \begin{defn}\label{def-return5}
Two minimal matchbox manifolds $\fM_i$ for $i=1,2$, are \emph{return equivalent} if there exist regular coverings of $\fM_i$ and  non--empty, clopen subsets $U_i \subset \fX_i$   so that the restricted pseudogroups $\cG^i_{U_i}$ for $i=1,2$ are return equivalent.
\end{defn}

This notion is well-defined by the following two results, whose proofs are given in Section~4 of \cite{CHL2013c}.

\begin{lemma}\label{lem-return1}
Let $\fM$ be a minimal matchbox manifold with transversal $\fX$ associated to a regular foliated covering of $\fM$.
Let  $W_1, W_2 \subset \fX$ be  non--empty open subsets,  then  there exist clopen subsets $U_i \subset W_i$ for $i=1,2$, such that 
the restricted pseudogroups $\cG^*_{U_1}$ and  $\cG^*_{U_2}$ are return equivalent.
\end{lemma}

\begin{lemma}\label{lem-return2}
Let $\fM$ be a minimal matchbox manifold, and suppose we are given  regular coverings by regular foliation charts 
 $\ds \{\vp_i \colon  U_i \to (-1,1)^n \times \fT_i \mid 1 \leq i \leq \nu\}$ and
 $\ds \{\vp_j' \colon  U_j' \to (-1,1)^n \times \fT_j' \mid 1 \leq j \leq \nu'\}$ of $\fM$, with transversals $\fX$ and $\fX'$ respectively.
 Let  $W_1  \subset \fX$ and $W_2  \subset \fX'$ be  non--empty open subsets. Then there exist clopen subsets $U_i \subset W_i$ for $i=1,2$, such that 
the restricted pseudogroups $\cG^*_{U_1}$ and  $\cG^*_{U_2}$ are return equivalent.
\end{lemma}

Finally, we have the   that homeomorphism implies return equivalence.
\begin{thm}\cite[Theorem~4.8]{CHL2013c}\label{thm-topinv}
Let $\fM_1$ and $\fM_2$ be minimal matchbox manifolds. Suppose that there exists a homeomorphism $h \colon \fM_1 \to \fM_2$, then $\fM_1$ and $\fM_2$  are return equivalent.
\end{thm}
 The question of when   return equivalence implies homeomorphism is also addressed in the work \cite{CHL2013c}.

We conclude this section with a fundamental observation about the dynamics of the pseudogroups associated to equicontinuous matchbox manifolds.
\begin{prop}\label{prop-inducedgroup}
Let   $\Phi \colon \cG \times \fX \to \fX$ be an equicontinuous pseudogroup action on a Cantor space $\fX$. 
Let  $W \subset \fX$ be a non--empty open subset,  then for any $x \in W$, there exists a clopen neighborhood $x \in U \subset W$ such  that the restricted $\psg$  $\cG^*_{U}$ is induced by a group action   $\Phi_U \colon \cH_U \times U \to U$, where  $\cH_U   \subset Homeo(U)$.
\end{prop}
\begin{remark}\label{rmk-codingpartition}
{\rm  
The conclusion of this result, as proved  in \cite[Section~6]{ClarkHurder2013}, actually includes a stronger statement. The   clopen set $U \subset \fX$ in Proposition \ref{thm-topinv} is constructed using the method of ``orbit coding'' for the pseudogroup action. Lemma~6.5 in   \cite{ClarkHurder2013} yields  that if $x \in U$ and $g \in \cG^*$ satisfies $\Phi(g)(x) \in U$ then there exists $h_{\cI} \in \cG^*$ defined as in \eqref{eq-defholo} such that $U \subset \fD(h_{\cI}) $ with $\Phi(h_{\cI})(U) = U$, and    $g = h_{\cI} | U \in \cH_U$. Every $g \in \cH_U$ is realized in this way. This implies, in particular, that the translates of $U$ under the action of $\cG^*$ are either contained in $U$, or disjoint from $U$. Thus,  the translates of $U$ under the action of $\cG^*$ form a clopen disjoint partition of $\fX$. The precise version of these facts is stated as Proposition~3.4 in \cite{DHL2016c}.}
\end{remark}

The action  $(\cH_U, U, \Phi_U)$ obtained in Proposition~\ref{prop-inducedgroup}  is an  \emph{equicontinuous minimal Cantor system} in the terminology of \cite{DHL2016a}.  

\begin{cor}\label{cor-inducedgroup}
Let $\fM_1$ and $\fM_2$ be equicontinuous matchbox manifolds, with transversals $\fX_1$ and $\fX_2$ defined by a choice of regular foliated coverings, respectively. Assume that their holonomy pseudogroups $\cG^*_{\F_1}$ and $\cG^*_{\F_2}$ are return equivalent. Then  there exist clopen subsets $U_i \subset \fX_i$   such that the restricted pseudogroup  $\cG^*_{U_i}$ is induced by a group action   $\Phi_{U_i} \colon \cH_{U_i} \times U_i \to U_i$, and 
  there exists a homeomorphism $h \colon U_1 \to U_2$ which induces an isomorphism of the   groups $\cH^1_{U_1}$ and $\cH^2_{U_2}$ and a
       conjugacy of equicontinuous minimal Cantor systems 
$(\cH^1_{U_1}, U_1, \Phi^1_{U_1})$ and $(\cH^2_{U_2}, U_2, \Phi^2_{U_2})$. 
\end{cor}

Clopen subsets with restricted pseudogroup action induced by a group action play an prominent role in the proofs throughout the paper, motivating the introduction of the following terminology.

\begin{defn}\label{defn-adapted}
Let   $\Phi \colon \cG \times \fX \to \fX$ be an equicontinuous pseudogroup action on a Cantor space $\fX$. 
We say that such a clopen subset $U \subset \fX$ is \emph{adapted} to the action $\Phi$, with group $\cH_U \subset Homeo(U)$, 
if  the restricted pseudogroup $\cG^*_U$ is induced by the group action $\Phi_U \colon  \cH_U \times U \to U$. 
\end{defn}

In the definition of return equivalence, there is no a priori restriction on the diameter of the clopen set  $U$, so it is important to investigate properties of the associated group  $\cH_{U}$ which are independent of the choice of the adapted clopen subset $U$. For example, the condition that $\cH_{U}$ is a virtually nilpotent group is such a property, which is a key idea in the work \cite{DHL2016b}.

 \section{Ellis group and the stable actions}\label{sec-ellis}

The  \emph{Ellis (enveloping) semigroup} associated to a continuous group action $\Phi \colon G \times X \to X$   was introduced in the papers \cite{EllisGottschalk1960,Ellis1960}, and is treated in the books   \cite{Auslander1988,Ellis1969,Ellis2014}.  We briefly recall some basic facts    for the special case of equicontinuous   minimal  Cantor systems, and then apply these concepts to the case of equicontinuous pseudogroup actions.  We  also recall   the notion of a \emph{strongly quasi-analytic} action introduced in the work \cite{ALC2009}, and then introduce the notions of \emph{stable} and \emph{wild} pseudogroup actions.

 \subsection{The Ellis group}\label{subsec-ellis}
  
 Let $X$ be a compact Hausdorff  topological space,   $G$ be a finitely generated group, and $\Phi \colon G \times X \to X$ an equicontinuous  action. 
 Let $\Phi(G) \subset Homeo(X)$ denote the image subgroup, and consider its closure $\overline{\Phi(G)}   \subset {\it Homeo}(X)$ in the \emph{uniform topology on maps}. 
 The group $\overline{\Phi(G)}$  is a special case of the more general construction of Ellis \emph{semi-groups} for topological actions.  Note that   each     $\whg \in \overline{\Phi(G)}$ is the limit of a sequence of maps in $\Phi(G)$; that is, we do not need to use ultrafilters in the definition of $\overline{\Phi(G)}$, in contrast to the more general construction of Ellis semi-groups.  
 By a small abuse of notation, we   use  the notation $\whg = (g_i)$ to denote a sequence $\{g_i \mid i \geq 1\} \subset G$ such that the sequence $\{\Phi(g_i) \mid i \geq 1\} \subset {\it Homeo}(X)$ converges to $\whg$.

If the action $\Phi \colon G \times X \to X$ is equicontinuous and minimal,   then  $\overline{\Phi(G)}(x) = X$ for any $x \in X$. That is, the group $\overline{\Phi(G)}$ acts transitively on $X$. Introduce the isotropy group  at $x$,  
\begin{align}\label{iso-defn2}
\overline{\Phi(G)}_x = \{ \whg=(g_i) \in \overline{\Phi(G)} \mid \whg( x) = x\},
\end{align}
then there is a   natural identification $X \cong \overline{\Phi(G)}/\overline{\Phi(G)}_x$ of left $G$-spaces.

The Ellis construction can be applied to the adapted clopen sets for an equicontinuous pseudogroup action $\Phi \colon \cG \times \fX \to \fX$ on a Cantor space $\fX$. 

 For   $x \in \fX$, let $U$ with $x \in U$ be a clopen set adapted to the action $\Phi$, with group  $\cH_{U}  \subset Homeo(U)$.
 Let  $\overline{\cH_U} \subset Homeo(U)$ denote the closure of $\cH_U$ in the uniform topology. Then, since the   action of $\cH_U$ on $U$ is minimal,  the action of $\overline{\cH_U}$ on $U$ is transitive, and so there is an $\overline{\cH_U}$-equivariant homeomorphism $U \cong \overline{\cH_U}/\overline{\cH_U}_x$ where $\overline{\cH_U}_x$ is the isotropy subgroup of  $x \in U$.
  We have the basic observation.
  \begin{lemma}\label{lem-closure-returnequiv}
For $i=1,2$,  let $\fU_i$   be a Cantor space,   let $H_i \subset Homeo(\fU_i)$ be a subgroup whose action is equicontinuous, and let $\overline{H_i} \subset Homeo(\fU_i)$ denote the closure of $H_i$ in the uniform topology on maps.
Suppose there exists a homeomorphism $\phi \colon  \fU_1 \to \fU_2$ which induces an isomorphism $\Phi_{\phi}: H_1 \to H_2$. 
 Then $\phi$     induces a topological isomorphism       $\ds \overline{\Phi}_{\phi}: \overline{H_1}  \to \overline{H_2}$. 
 \end{lemma}

 \subsection{Stable and wild actions}\label{subsec-defSQA}
 
The  \emph{strong quasi-analyticity} condition for a pseudogroup action   was introduced by  A.~{\'A}lvarez L{\'o}pez and A.~Candel  in their work \cite{ALC2009}, and   has a fundamental role in their construction of Molino spaces  for equicontinuous foliated spaces in \cite{ALM2016} (see  \cite[Definition~2.18]{ALM2016} and the discussion of the concept there). Their definition, stated below,  was motivated by the search for a condition equivalent to quasi-analyticity condition,  as introduced by Haefliger \cite{Haefliger1985}  for the $\psg$s of  smooth foliations.

 \begin{defn}\cite{ALC2009} \label{def-sqa1}
 A $\psg$ $\cG^*$ acting on a locally compact space $\fX$ is \emph{strongly quasi-analytic}, or \emph{SQA}, if for every $h \in \cG^*$ the following holds: Let $U \subset \fD(h)$ be an open set, and suppose that the restriction $h|U$ is the identity map, then $h$ is the identity map on   $\fD(h)$.
 \end{defn}
 Note that an analytic pseudogroup, as defined by Definition~\ref{def-analytic}, satisfies the SQA condition, though an SQA action need not be analytic. 
A group action $\Phi \colon G \times \fX \to \fX$ is \emph{free} if for any $x \in \fX$ and $g \in G$, then $\Phi(g)(x) = x$ implies that $g$ is the identity element of $G$. It is   immediate from the definitions that a free action satisfies the SQA condition.

 For the study of matchbox manifolds,    there is no preferred choice of a transverse section $\fX$, and so the following \emph{local} SQA property is a more natural property to require. 

\begin{defn}\label{def-localsqa}
 A $\psg$ $\cG^*$ acting on a locally compact space $\fX$ is \emph{locally strongly quasi-analytic}, or \emph{LSQA}, if for each $x \in \fX$, there exists an open subset $U \subset \fX$ with $x \in U$ such that the restricted $\psg$ $\cG^*_U$, as defined by \eqref{eq-inducedpseudo}, is strongly quasi-analytic.
\end{defn}

It is possible that the restricted holonomy $\psg$ of a weak solenoid is SQA for a suitable choice of a   section, though the global holonomy action of the solenoid is not SQA. This is true, for example,   in the case of  \cite[Example 7.6]{DHL2016a}.  Thus, the global holonomy action of a weak solenoid may be LSQA, though need not satisfy the SQA condition for a given section $\fX$.

 A basic issue with the SQA and the LSQA properties, as formulated in Definitions~\ref{def-sqa1} and \ref{def-localsqa}, is that for the applications of these properties to the Molino theory for equicontinuous foliated spaces, one also needs this property for the closure of a  pseudogroup action in the compact-open topology of maps  for their restriction to open subsets of the transversal.   We formulate this as follows.
 
 \begin{defn}\label{def-LCSQA}
 A $\psg$ $\cG^*$ acting on a locally compact space $\fX$ is \emph{locally completely strongly quasi-analytic}, or \emph{LCSQA}, if for each $x \in \fX$,    there exists an   open subset with $x \in U$  such that  for any $\whg \in \overline{\cG^*_U}$, if there exists an open set $V \subset Dom(\whg) \cap U$ for which  the restriction of $\whg$ to $V$ is the identity map, then $\whg$ is the identity map on its domain.
 \end{defn}

 The   notion of an LCSQA pseudogroup action has a   more straightforward formulation in the case where   $\fX$ is a Cantor space, and so its topology is generated by clopen sets. This yields the notion of a \emph{stable} equicontinuous pseudogroup action $\Phi \colon \cG \times \fX \to \fX$  on a Cantor space $\fX$ as follows.
 
Consider a clopen set $U \subset \fX$ adapted to the action $\Phi$, with group $\cH_{U}$. Then for $x \in U$, let $V \subset U$ be a subset, adapted to the action $\Phi$, with group $\cH_{V} $. 
Let $\cH_{U,V} \subset \cH_U \subset Homeo(U)$ denote the subgroup of elements $g \in \cH_U$ such that $g ( V )= V$.
Since $V$ is clopen in $U$, the closure $\overline{\cH_{U,V}}$ consists of   elements $\whg \in \overline{\cH_U}$ such that $\whg ( V )= V$.
Then there exists   homomorphisms induced by restriction,
\begin{equation}\label{eq-restriction}
\rho_{U,V} \colon \cH_{U,V}\to \cH_{V} \quad , \quad \orho_{U,V} \colon \overline{\cH_{U,V}} \to \overline{\cH_{V}} \ .
\end{equation}
  
   \begin{defn} \label{def-stable2}
Let   $\Phi \colon \cG \times \fX \to \fX$ be an equicontinuous pseudogroup action on a Cantor space $\fX$. We say that the action is \emph{stable}  
if for each $x \in \fX$, there exists a clopen set $U \subset \fX$ adapted to the action $\Phi$,  with $x \in U$ and with group $\cH_{U} $, such  that for any
clopen set $V \subset U$ adapted to $\Phi$ with group  $\cH_{V}$, 
the restriction map $\orho_{U,V}$ \emph{on closures}   in \eqref{eq-restriction} has trivial kernel.

The action is said to be \emph{wild} if it is not stable.
 \end{defn}

 We use Proposition~\ref{prop-inducedgroup}  to show the following:

 \begin{prop}\label{prop-sqaequiv}
 An equicontinuous pseudogroup action $\Phi \colon \cG \times \fX \to \fX$ on a Cantor space $\fX$ is stable if and only if it is LCSQA.
  \end{prop}
 
 \proof Suppose the action $\Phi$ is stable. Given $x \in \fX$, by Proposition~\ref{prop-inducedgroup}, there is a clopen subset  $x \in U \subset \fX$ adapted to $\Phi$ with group $\cH_U$. We show that the action of  the closure $\overline{\cH_U}$ on $U$ is SQA.

Let $\whg \in  \overline{\cH_U}$ be  such that for some open subset $V \subset U$,  the restriction $\whg|V$ is the identity map.  Again by Proposition~\ref{prop-inducedgroup}, there exists a non-empty clopen subset $W \subset V$, adapted to $\Phi$, with group $\cH_{W} $. 
 By the assumption that the action $\Phi$ is stable, the map $\orho_{U,W}$ has trivial kernel, and as $\whg$ restricts to the identity map on $W$, it follows that $\whg$ is the identity map on $U$, as was to be shown.

 Now suppose the action  of $\cG^*$ on $\fX$ is LCSQA. Let  $x \in \fX$, then there exists a  clopen neighborhood $x \in W \subset \fX$ such that   the   restricted pseudogroup $\cG^*_W$ satisfies the SQA condition for $\overline{\cG^*_W}$. 
Then  by Proposition~\ref{prop-inducedgroup}, there exists a non-empty clopen subset $U \subset W$ with $x \in U$, adapted to the action $\Phi$, with group $ \cH_{U} $.   

Let $V \subset U$ be a clopen set, adapted to $\Phi$. We   show that the kernel of the restriction map $\orho_{U,V} \colon \overline{\cH_{U,V}} \to \overline{\cH_{V}}$   is trivial.  
 Let $\whg \in \overline{\cH_{U,V}}$ be such that its restriction to $V$ is the identity map. Then as $\cG^*_W$ satisfies the LCSQA condition, $\whg$ must be the identity on   $U$.
  \endproof

The definition of a stable action given in \cite{DHL2016c} was formulated in terms of group chains, and we show in Proposition~\ref{prop-stableisstable} below that the definition in  \cite{DHL2016c}  coincides with the notion in Definition~\ref{def-stable2}. It was shown in \cite{DHL2016c}, that the holonomy pseudogroup associated to the Schori solenoid is not LSQA, and so is not LCSQA.  As a consequence of    Proposition~\ref{prop-sqaequiv}, we obtain that the Schori weak solenoid as described   in Theorem~9.8 of \cite{DHL2016c} is an example of a non-stable equicontinuous minimal group action.
 
\begin{cor}\label{cor-Schori-wild}
The global monodromy action of the Schori solenoid  is wild.
\end{cor}

Finally, we show  that the stable   property is preserved by return equivalence.   The following result is the analog of Lemma~9.5 in \cite{ALC2009}.

\begin{prop}\label{prop-retequivstable}
The stable property is an invariant of return equivalence for   equicontinuous minimal pseudogroup Cantor actions.
\end{prop}
\proof
 For $i =1,2$, let   $\Phi^i \colon \cG_i \times \fX_i \to \fX_i$  be a minimal  equicontinuous pseudogroup action  on  a Cantor space $\fX_i$, and 
 assume that the actions are return equivalent. Assume that the action of $\cG_1$ on $\fX_1$ is stable, then we show that the action of $\cG_2$ on $\fX_2$ is also stable.

 For   $i=1,2$, by the Definition~\ref{def-return},  there exist clopen subsets $W_i \subset \fX_i$, adapted to $\Phi_i$,  and a homeomorphism $\phi \colon W_1 \to W_2$ which induces an isomorphism $\Phi_{\phi} \colon \cH^1_{W_1} \to \cH^2_{W_2}$. Then by Lemma~\ref{lem-closure-returnequiv}, there is an induced isomorphism between their closures, 
 $\ds \overline{\Phi}_\phi \colon  \overline{\cH^1_{W_1} } \to \overline{\cH^2_{W_2}}$. 
 
 By Proposition~\ref{prop-inducedgroup}, we can choose an adapted clopen subset $U_1 \subset W_1$  such that $x \in U_1$, with group $\cH_{U_1}^1$. 
The image $U_2 = \phi(U_1)$ is a clopen subset of $\fX_2$. 
By Lemma~\ref{lem-finitegen},  there exists an integer  $m$ such that $\fX_2$ is covered by the finite collection
$\cU_2 = \{h_{\cI_1}(U_2), \ldots, h_{\cI_m}(U_2)\}$ of clopen sets in $\fX_2$ where each $h_{\cI_{\ell}} \in \cG^*_2$. 

Given $y \in \fX_2$ there exists an index $1 \leq i_y \leq m$ such that $y \in h_{\cI_{i_y}}(U_2)$. 
By Proposition~\ref{prop-inducedgroup}, we can choose a clopen set $U_y$, adapted to $\Phi_2$, such that  $y \in U_y \subset h_{\cI_{i_y}}(U_2)$, with group   $\cH^2_{U_y}$. 

  Note that $h_{{\cI_{i_y}}}^{-1}(U_y) \subset U_2$, and we set $U_y' =   \phi^{-1} ( h_{{\cI_{i_y}}}^{-1}(U_y)) \subset U_1$.
Then $h_{\cI_{\ell}} \circ \phi$ induces a topological conjugacy between $\cH^1_{U_y'}$ and $\cH^2_{U_y}$,  and by Lemma~\ref{lem-closure-returnequiv}  between their closures in the uniform topologies.

Now let $V \subset U_y$ be a clopen subset, adapted to $\Phi_2$, with group $ \cH^2_{V} $. 
Then $V' =   \phi^{-1} ( h_{{\cI_{i_y}}}^{-1}(V)) \subset U_y' \subset U_1$ which is again a clopen subset.

Let $\whg \in \overline{\cH^2(U_y , V)}$ be such that its restriction to $V$ is the identity map, then it is conjugate by  $h_{\cI_{\ell}} \circ \phi$   to a map   $\whg' \in \overline{\cH^1_{U_y', V'}}$ which restricts to the identity map on $V'$.  As the action of $\cG_1$ on $\fX_1$ is stable, $\whg'$ is the identity on $U_y'$, hence $\whg$ is the identity on $U_y$. Thus,   the action of $\cG_2$ on $\fX_2$ is   stable.

 Reversing the roles of $\cG_1$ and $\cG_2$   shows the converse implication, and completes the proof.
 \endproof

 \eject
 
\section{The asymptotic discriminant invariant}\label{sec-adi}

 In this   section, we introduce the notion of tail equivalence for a sequence of group homomorphisms, and apply the results of the previous section   to obtain a new invariant of the homeomorphism type for equicontinuous matchbox manifolds.
 The results of this section can be seen as an extension of those of the last section, and several of the results are analogous to results in Section~\ref{sec-ellis}. 
We begin with a definition, which is adapted from the notion of equivalence for chains of covering spaces as defined by Rogers and Tollefson in \cite[Section 3]{RT1971b}, and formulated for group chains  by Fokkink and Oversteegen in \cite[Definition~14]{FO2002}.

  \begin{defn} \label{def-tailequiv}
 Let $\cA = \{\phi_{i} \colon A_{i} \to A_{i+1} \mid i \geq 1\}$ and $\cB = \{\psi_{i} \colon B_{i} \to B_{i+1} \mid i \geq 1\}$ be two  sequences of  \emph{surjective}   group homomorphisms. 
 Say that $\cA$ and $\cB$ are tail equivalent, and write $\cA \tail \cB$, if the sequences of groups  $\cA$ and $\cB$ are intertwined 
 by a   sequence of surjective   group homomorphisms. 
 That is, there exists:
 \begin{enumerate}
\item an increasing sequence of indices $\{p_i   \mid i \geq   1 ~ , ~ p_{i+1} > p_{i} \geq i \geq 1 \}$;
\item an increasing sequence of indices $\{q_i   \mid i \geq   1 ~ , ~ q_{i+1} > q_{i} \geq i \geq 1 \}$;
\item a sequence $\cC = \{ \tau_i \colon C_i \to C_{i+1} \mid i \geq 1\}$ of surjective group homomorphisms;
\item a collection of   isomorphisms $\Pi_{\cA\cC} \equiv \{ \Pi^{i}_{\cA \cC} \colon A_{p_i} \to C_{2i-1}\mid i \geq  1\}$;
\item a collection of   isomorphisms $\Pi_{\cB\cC} \equiv \{ \Pi^{i}_{\cB \cC} \colon B_{q_i} \to C_{2i}\mid i \geq  1\}$;
\end{enumerate}
such that  for all $i \geq 1$, we have 
\begin{eqnarray} 
 \tau_{2i} \circ   \tau_{2i-1} \circ \Pi^i_{\cA\cC}   & = &  ~    \Pi^{i+1}_{\cA\cC}  \circ \Phi_{i}  \ ,   \label{eq-interlacedAB1} \\ 
 \tau_{2i+1} \circ   \tau_{2i} \circ \Pi^i_{\cB\cC}   & = &  ~    \Pi^{i+1}_{\cB\cC}  \circ \Psi_{i}  \ ,  \label{eq-interlacedAB2}
 \end{eqnarray} 
 where 
 \begin{eqnarray}
\Phi_{i} & = &  \phi_{p_{i+1} -1} \circ  \phi_{p_{i+ 1} -2} \circ \cdots \circ \phi_{p_i +1 } \circ  \phi_{p_i} \\
\Psi_{i} & = &   \psi_{q_{i+1} -1} \circ  \psi_{q_{i+ 1} -2} \circ \cdots \circ \psi_{q_i +1 } \circ  \psi_{q_i}   \ .
\end{eqnarray}
 \end{defn}
The collections of maps $\Pi_{\cA\cC}$  and $\Pi_{\cB\cC}$ satisfying \eqref{eq-interlacedAB1} and \eqref{eq-interlacedAB2} are said to  \emph{realize}  $\cA \tail \cB$. 
 The relations \eqref{eq-interlacedAB1} and \eqref{eq-interlacedAB2} give rise to a commutative diagram of maps

    \begin{align} \label{eq-commutativediagram}
 \xymatrix{
 \cdots A_{p_i}   \ar[r]^{\phi_{p_i}} \ar[d]_{\Pi^i_{\cA\cC}}^{\cong} &\cdots  \ar[r]^{\phi_{p_{i+1} -1}} &    A_{p_{i+1}} \ar[r]^{\phi_{p_{i+1}}} \ar[d]_{\Pi^{i+1}_{\cA\cC}}^{\cong} & \cdots  \ar[r]^{\phi_{p_{i+2} -1}}&   A_{p_{i+2}} \ar[r]^{\phi_{p_{i+2}}}  \ar[d]_{\Pi^{i+2}_{\cA\cC}}^{\cong} &     \cdots \\
\cdots C_{2i-1}   \ar[r]_{\tau_{2i-1}}    & C_{2i}   \ar[r]_{\tau_{2i}}  &          C_{2i+1} \ar[r]_{\tau_{2i+1}}   & C_{2i+2}   \ar[r]_{\tau_{2i+2}} & C_{2i+3}    \ar[r]_{\tau_{2i+3}}  &  C_{2i+4}  & \cdots      \\
\cdots     \ar[r]_{\psi_{q_i -1}}    & B_{q_{i }}   \ar[r]_{\psi_{q_i }} \ar[u]_{\Pi^{i}_{\cB\cC}}^{\cong} &    \cdots  \ar[r]_{\psi_{q_{i+1} -1}}   &      B_{q_{i+1}} \ar[r]_{\psi_{q_{i+1}}}   \ar[u]_{\Pi^{i+1}_{\cB\cC}}^{\cong} & \cdots  \ar[r]_{\psi_{q_{i+2} -1}}  & B_{q_{i+2}}  \ar[u]_{\Pi^{i+2}_{\cB\cC}}^{\cong} &  \cdots 
} 
\end{align}

 We  omit the standard proof that ``tail equivalence''  is an equivalence relation. The tail equivalence class of a sequence $\cA$ is denoted by $[\cA]_{\infty}$ and is called the \emph{asymptotic class} of $\cA$. 
 
 We also have the usual notion of isomorphisms of sequences.
  \begin{defn} \label{def-isoseq}
  Let $\cA = \{\phi_{i} \colon A_{i} \to A_{i+1} \mid i \geq 1\}$ and $\cB = \{\psi_{i} \colon B_{i} \to B_{i+1} \mid i \geq 1\}$ be   sequences of  surjective   group homomorphisms. 
 Say that $\cA$ and $\cB$ are isomorphic, and write $\cA \cong \cB$, if there exists isomorphisms $h_i \colon A_i \to B_i$ such that 
 $\psi_i \circ h_i = h_{i+1} \circ \phi_i$ for all $i \geq 1$.
 \end{defn}

A sequence $\cA$ is \emph{constant} if   each map $\phi_{i} \colon A_{i} \to A_{i+1}$ is an isomorphism for all $i \geq 1$, and $\cB$ said to be \emph{asymptotically constant} if it is tail equivalent to a constant sequence $\cA$. The following  result follows from the usual method of ``chasing of diagrams''.

\begin{lemma}\label{lem-asymptoticconstant}
A sequence $\ds \cB = \{\psi_{i} \colon B_{i} \to B_{i+1} \mid i \geq 1\}$ of  surjective homomorphisms is asymptotically constant if and only if there exists $i_0 \geq 0$ such that $ker(\psi_{i})$ is trivial for all $i \geq i_0$.
\end{lemma}

We next apply the notion of tail equivalence to obtain a new invariant for the return equivalence class of  an equicontinuous minimal Cantor pseudogroup action 
  $\Phi \colon \cG \times \fX \to \fX$ on a Cantor space $\fX$. Let $\cG^*$ denote the associated $\psg$ acting on $\fX$.
  
  Recall that by Definition \ref{defn-adapted}, a clopen subset  $U \subset \fX$, $x \in U$, is \emph{adapted} to the action $\Phi$  if 
the restricted pseudogroup   $\cG^*_U$ is induced by a group action $\Phi_U \colon  \cH_U \times U \to U$, with image group $\cH_U \subset Homeo(U)$. Recall from Section \ref{subsec-ellis} that the closure $\overline{\cH_U} \subset Homeo(U)$ of $\cH_U$ in the uniform topology on maps is isomorphic to the Ellis group of the group action $\Phi_U$. Also, since $\Phi_U$ is minimal, the closure $\overline{\cH_U} $   acts transitively on $U$.

We will need the following basic observation about adapted sets and their properties.

\begin{lemma}\label{lem-adaptedsurjective}
Let  $\Phi \colon \cG \times \fX \to \fX$ be an equicontinuous minimal pseudogroup action on a Cantor space $\fX$.
Let  $V \subset U \subset \fX$  be    clopen sets, and assume that both $U$ and $V$ are adapted to the action $\Phi$.  Then  the restriction map
$\ds \orho_{U,V} \colon \overline{\cH_{U,V}} \to \overline{\cH_{V}}$ defined in \eqref{eq-restriction} is a surjection. 
\end{lemma}
\proof
Let $\whh \in \overline{\cH_{V}}$ then there exists a sequence $\{h_{i} \in \cH_{V} \mid i \geq 1\}$ with $h_{i} \to \whh$ in the uniform topology on maps of $V$. As $V \subset U$ is an open set, the maps $h_i$ are elements in $\cG^*_{U}$. As $\cG^*_{U}$ is induced from the group action    $\Phi_{U} \colon  \cH_{U} \times U  \to U$  there exists some $g_i \in \cH_{U}$ whose restriction to $V$ satisfies 
$g_i | V = h_i$. The group $\overline{\cH_{U}}$ is compact, so there exists a subsequence $g_{i_j}$ which converges to $\whg \in \overline{\cH_{U}}$, and also the subsequence $g_{i_j} | V$  converges to $\whh$. Thus, $\whg | V = \whh$ as was to be shown.
\endproof

\begin{defn}\label{def-adaptednbhds}
Let  $\Phi \colon \cG \times \fX \to \fX$ be an equicontinuous minimal pseudogroup action on a Cantor space $\fX$.
Then a descending chain of clopen sets $\cU = \{U_{\ell} \subset \fX  \mid \ell \geq 1\}$ is said to be an \emph{adapted neighborhood basis} at $x \in \fX$ for the action $\Phi$  if
    $x \in U_{\ell +1} \subset U_{\ell}$ for all $ \ell \geq 1$ with     $\cap  \ U_{\ell} = \{x\}$, and  each $U_{\ell}$ is adapted to the action $\Phi$.
\end{defn}

Proposition~\ref{prop-inducedgroup} implies that each $x \in \fX$ has an adapted   neighborhood basis.

\begin{lemma}\label{lem-restrictions}
Let  $\Phi \colon \cG \times \fX \to \fX$ be an equicontinuous minimal action on a Cantor space $\fX$, and 
let $\cU = \{U_{\ell} \subset \fX  \mid \ell \geq 1\}$ be an \emph{adapted neighborhood basis} at $x \in \fX$. Let $\ell' > \ell > 1$, so that $U_{\ell'} \subset U_{\ell}$. Suppose that  $\whh \in \overline{\cH_{U_1}}$ and  $\whh(U_{\ell'}) = U_{\ell'}$,  then    $\whh(U_{\ell}) = U_{\ell}$.
\end{lemma}
\proof
By Remark~\ref{rmk-codingpartition}, the translates of $U_{\ell}$ under the action of $\cH_{U_1}$ form a disjoint clopen partition of $U_1$. Suppose that   for $g \in G$ we have $\Phi(g)(U_{\ell'}) = U_{\ell'}$ then $\Phi(g)(U_{\ell}) \cap U_{\ell} \ne \emptyset$, hence $\Phi(g)(U_{\ell}) = U_{\ell}$. 
 
Given $\whh \in \overline{\cH_{U_{1}}}$ there exists a sequence $\{g_{i} \mid i \geq 1\} \subset G$ such that the maps $\Phi(g_{i}) | U_{1}$ converge in the uniform topology to $\whh$.
As each set in the disjoint partition of $U_1$ by the translates of $U_{\ell'}$ is clopen, for   $i$ sufficiently large we must have $\Phi(g_{i})(U_{\ell'}) = U_{\ell'}$, 
and hence $\Phi(g_{i})(U_{\ell}) = U_{\ell}$. As $U_{\ell}$ is compact,  we conclude that   $\whh(U_{\ell}) = U_{\ell}$.
\endproof

 Let  $\Phi \colon \cG \times \fX \to \fX$ be an equicontinuous   pseudogroup action on a Cantor space $\fX$.  
Consider an adapted neighborhood system $\cU = \{U_{\ell} \subset \fX  \mid \ell \geq 1\}$ at $x \in \fX$.
For $\ell' > \ell \geq 1$, following the constructions  in Section~\ref{subsec-defSQA}, define  $\cH_{U_{\ell},U_{\ell'}} \subset \cH_{U_{\ell}} \subset Homeo(U_{\ell})$ to be the subgroup of elements $h \in \cH_{U_{\ell}}$ such that $h (U_{\ell'})= U_{\ell'}$, with closure   
$\overline{\cH_{U_{\ell},U_{\ell'}}} \subset  \overline{\cH_{U_{\ell}}}$. Then for $\whh \in \overline{\cH_{U_{\ell},U_{\ell'}}}$ we have $\whh(U_{\ell'}) = U_{\ell'}$, so there are well-defined  restriction homomorphisms
\begin{equation}\label{eq-restriction33}
\rho_{U_{\ell} , U_{\ell'}} \colon \cH_{U_{\ell},U_{\ell'}}\to \cH_{U_{\ell'}} \quad , \quad \orho_{U_{\ell} , U_{\ell'}} \colon \overline{\cH_{U_{\ell} , U_{\ell'}}} \to \overline{\cH_{U_{\ell'}}} \ .
\end{equation}

 Introduce the isotropy groups of the actions of the closure $ \overline{\cH_{U_1 , U_{\ell}}}$, given by
 \begin{equation}\label{eq-defisotropy}
\cI(U_1, U_{\ell},x)  = \overline{\cH_{U_1 , U_{\ell}}}_x = \{ \whh \in \overline{\cH_{U_1 , U_{\ell}}} \mid \whh(x) = x\} \ .
\end{equation}
Then we have the following consequence of Lemma~\ref{lem-restrictions}.
\begin{cor}\label{cor-restrictionsrho}
Let  $\Phi \colon \cG \times \fX \to \fX$ be an equicontinuous   pseudogroup action on a Cantor space $\fX$, and 
let $\cU = \{U_{\ell} \subset \fX  \mid \ell \geq 1\}$ be an  adapted neighborhood basis  at $x \in \fX$.
For $\ell' > \ell \geq 1$,   the restriction map
$\ds \orho_{U_{\ell},U_{\ell'}}$ induces a surjection $\sigma_{U_{\ell},U_{\ell'}} \colon \cI(U_1, U_{\ell},x) \to \cI(U_1, U_{\ell'},x)$.
\end{cor}
 \proof
 Let $\whh \in \cI(U_1, U_{\ell'},x) \subset  \overline{\cH_{U_{1}}}$ then $\whh(U_{\ell'}) = U_{\ell'}$. Then by Lemma~\ref{lem-restrictions} we have 
 $\whh(U_{\ell}) = U_{\ell}$ so that $\whh \in \cI(U_1, U_{\ell},x)$ and thus $\sigma_{U_{\ell},U_{\ell'}}$ is onto. 
 \endproof

 Let  $\Phi \colon \cG \times \fX \to \fX$ be an equicontinuous minimal pseudogroup action on a Cantor space $\fX$, let $x \in \fX$, and 
let $\cU = \{U_{\ell} \subset \fX  \mid \ell \geq 1\}$ be an  adapted neighborhood basis  at $x$.
By Corollary~\ref{eq-defisotropy}, we have  the sequence of surjective homomorphisms, 
\begin{equation}\label{eq-surjectiveUn}
\cI(\cU,x) = \{\sigma_{U_{\ell},U_{\ell+1}} \colon \cI(U_1, U_{\ell},x) \to \cI(U_1, U_{\ell+1},x) \mid \ell \geq 1\} \ .
\end{equation} 
The following results  consider the tail equivalence class $[\cI(\cU,x)]_{\infty}$.
 \begin{lemma}\label{lem-adaptedinv1}
  Let  $\Phi \colon \cG \times \fX \to \fX$ be an equicontinuous   pseudogroup action on a Cantor space $\fX$, and 
let $\cU = \{U_{\ell} \subset \fX  \mid \ell \geq 1\}$ and $\cV = \{V_{\ell} \subset \fX  \mid \ell \geq 1\}$ be adapted neighborhood bases  at $x \in \fX$.
Then the sequences of surjective homomorphisms $\cI(\cU,x)$ and $\cI(\cV,x)$ are tail equivalent. 
  \end{lemma}
 \proof
 Let $\cU = \{U_{\ell} \subset \fX  \mid \ell \geq 1\}$ and $\cV = \{V_{\ell} \subset \fX  \mid \ell \geq 1\}$ be adapted neighborhood bases  at $x$.
 We introduce a third neighborhood basis $\cW = \{W_{\ell} \subset \fX  \mid \ell \geq 1\}$, defined recursively as follows.
 
 Set $W_1 = U_1$ and $p_1 = 1$. Then there exists $q_1 \geq 1$ such that $V_{q_1} \subset U_1$ and we set $W_2 = V_{q_1}$. Now proceed recursively, and assume that $\{p_1 , p_2, \ldots , p_i\}$ and $\{q_1, q_2, \ldots , q_i\}$ have been chosen, where   $W_{2j-1} = U_{p_j}$ and $W_{2j} = V_{q_j}$ for $1 \leq j \leq i$. Then choose $p_{i+1} > p_i$  such that $U_{p_{i+1}} \subset V_{q_i}$ and set $W_{2j+1} = U_{p_{i+1}}$.  
 Choose $q_{i+1} > q_i$ such that $V_{q_{i+1}} \subset U_{p_{i+1}}$ and set $W_{2i+2} = V_{q_{i+1}}$. 
 
Next, let $A_{i} = \cI(U_1, U_{i},x)$, $B_{i} = \cI(V_1, V_{i},x)$, and $C_{i} = \cI(U_1, W_{i},x)$. 

Define $\phi_i = \sigma_{U_{i},U_{i+1}}$,  $\psi_i = \sigma_{V_{i},V_{i+1}}$, and $\tau_{i} = \sigma_{W_i,W_{i+1}}$.

Define  $\Pi^i_{\cA\cC} \colon A_{p_i} \to C_{2i-1} = \cI(U_1, W_{2i-1},x) = \cI(U_1, U_{p_i},x) = A_{p_i}$ to be the identity for each $i \geq 1$.

Define  $\Pi^i_{\cB\cC} \colon B_{q_i} \to C_{2i} = \cI(U_1, W_{2i},x) = \cI(U_1, V_{q_i},x) = B_{q_i}$ to be the identity for each $i \geq 1$.

Then the identities \eqref{eq-interlacedAB1} and \eqref{eq-interlacedAB2} as in the diagram \eqref{eq-commutativediagram} are satisfied, so that $\cI(\cU,x)$ and $\cI(\cV,x)$ are tail equivalent, as was to be shown. 
\endproof

\begin{lemma}\label{lem-adaptedinv2} 
Let $\Phi \colon \cG \times \fX \to \fX$ be an equicontinuous pseudogroup action on a Cantor space $\fX$. 
Let $x \in \fX$, and suppose that $\cU = \{U_{\ell} \subset \fX  \mid \ell \geq 1\}$ is an  adapted neighborhood basis  at $x \in \fX$.
Let $h \in \cG^*$ with $h(x) = y \in \fX$, and let $\cV = \{V_{\ell} \subset \fX  \mid \ell \geq 1\}$ be an  adapted neighborhood basis  at $y$.
Then the sequences of surjective homomorphisms $\cI(\cU,x)$ and $\cI(\cV,y)$ are tail equivalent. 
\end{lemma}
\proof 
Recall that $\fD(h) \subset \fX$ denotes the domain of $h$ which is an open neighborhood with $x \in \fD(h)$, so there exists $\ell_0 \geq 1$ such that $U_{\ell} \subset \fD(h)$ for $\ell \geq \ell_0$. Each image $h(U_{\ell}) \subset \fX$ is then a clopen neighborhood of $y$, so we obtain a 
   chain of clopen subgroups  $\ds \cU^h = \{h(U_{\ell})    \mid \ell \geq \ell_0\}$ at $y$. Moreover, as each $U_{\ell}$ is an adapted clopen set, the same holds for the image $h(U_{\ell})$. We then have that $\cI(\cU,x)$ and $\cI(\cU^h,y)$ are tail equivalent by construction, and $\cI(\cU^h,y)$ is tail equivalent to $\cI(\cV,y)$ by 
  Lemma~\ref{lem-adaptedinv1}. Thus,   $\cI(\cU,x)$ and $\cI(\cV,y)$  are tail equivalent, as was to be shown.
\endproof

 \begin{lemma}\label{lem-adaptedinv3}
 Let $\Phi \colon \cG \times \fX \to \fX$ be an equicontinuous \emph{minimal} pseudogroup action on a Cantor space $\fX$, and 
assume that $\cU = \{U_{\ell} \subset \fX  \mid \ell \geq 1\}$ is an  adapted neighborhood basis  at $x \in \fX$.
Let $y \in U_1$ and  $\cV = \{V_{\ell} \subset \fX  \mid \ell \geq 1\}$ be an adapted neighborhood basis at $y$. 
Then   the sequences of surjective homomorphisms $\cI(\cU,x)$ and $\cI(\cV,y)$ are tail equivalent. 
\end{lemma}
\proof
The action $\Phi$ is minimal, so the restricted action $\Phi|U_1$ is also minimal. As $U_1$ is an adapted clopen set, this implies that the action of the closure $\overline{\cH_{U_1}}$ is transitive. Thus, there exists $\whh \in \overline{\cH_{U_1}}$ such that $\whh(x) = y$. 
Introduce the sequence of clopen  subsets of  $U_1$ given by  $\cU^\whh = \{\whh(U_{\ell})    \mid \ell \geq  2\}$ which is a neighborhood basis at $y$.

We claim that   $\whh(U_{\ell})$ is adapted to the action $\Phi$, for each $\ell \geq 2$.  
Choose a sequence $\{h_i \in \cH_{U_1}\}$ such that $h_i$ converges to  $\whh$ in the uniform topology on maps. 

 Fix $\ell \geq 2$, then $U_{\ell}$ is a clopen set and $\whh \in Homeo(U_1)$ implies that $\whh(U_{\ell})$ is a clopen subset of $U_1$ 
 and hence $\e_{\ell} = \dX\{\whh(U_{\ell}) , U_1 - \whh(U_{\ell})\} > 0$. Thus, 
  there exists $i_{\ell} \geq 1$ such that $i \geq i_{\ell}$ 
implies that $h_i(U_{\ell}) \cap \fX - \whh(U_{\ell}) = \emptyset$, 
and hence $h_i(U_{\ell}) = \whh(U_{\ell})$. 
As $h_i \in \cH_{U_1}$ this implies that $\whh(U_{\ell})$ is adapted to the action of $\Phi$.
 Thus,  $\cU^\whh$ is an adapted neighborhood basis at $y$.
 
 By Lemma~\ref{lem-adaptedinv1} the sequences of surjective homomorphisms $\cI(\cU^{\whh},y)$ and $\cI(\cV,y)$ are tail equivalent. We   have by construction that $\cI(\cU,x)$  and $\cI(\cU^{\whh},y)$  are tail equivalent, so $\cI(\cU,x)$ is tail equivalent to  $\cI(\cV,y)$  as claimed.
 \endproof

  \begin{lemma}\label{lem-adaptedinv4}
 Let $\Phi \colon \cG \times \fX \to \fX$ be an equicontinuous minimal pseudogroup action on a Cantor space $\fX$.
 Let $\cU = \{U_{\ell} \subset \fX  \mid \ell \geq 1\}$  and 
 $\cV = \{V_{\ell} \subset \fX  \mid \ell \geq 1\}$ be   adapted neighborhood bases at $x \in \fX$ and $y\in \fX$, respectively.
 Then   the sequences of surjective homomorphisms $\cI(\cU,x)$ and $\cI(\cV,y)$ are tail equivalent.
 \end{lemma}
\proof
By Lemma~\ref{lem-finitegen} there exists $h  \in \cG^*$ such that $y \in h(U_1 \cap \fD(h))$. Let $z = h^{-1}(y) \in U_1 \cap \fD(h)$. 
Then by  Proposition~\ref{prop-inducedgroup}, there exists an  adapted neighborhood basis $\cW = \{W_{\ell} \subset U_1 \cap \fD(h)  \mid \ell \geq 1\}$ at $z$.
Then by Lemma~\ref{lem-adaptedinv3}, the sequences   $\cI(\cU ,x)$ and $\cI(\cW,z)$ are tail equivalent.

Let $\ell_y \geq 1$ be such that $h(W_{\ell}) \subset V_1$ for $\ell \geq \ell_y$. 
Define $\cW^h = \{h(W_{\ell}) \subset V_1 \mid \ell \geq \ell_y\}$. Then the sequences $\cI(\cW, z)$ and $\cI(\cW^h, y)$ are tail equivalent by construction.
Then by Lemma~\ref{lem-adaptedinv1}, the sequences $\cI(\cW^h, y)$ and $\cI(\cV, y)$ are tail equivalent.
 
 Thus, we have $\cI(\cU ,x) \tail \cI(\cW, z) \tail \cI(\cW^h, y) \tail \cI(\cV, y)$ as was to be shown.
\endproof

As a consequence of  Lemma~\ref{lem-adaptedinv4} we obtain a well-defined invariant of the action $\Phi$.

\begin{defn}\label{def-asympdisc}
   Let $\Phi \colon  \cG \times \fX \to \fX$ an equicontinuous minimal pseudogroup action  on a Cantor space $\fX$.
The   \emph{asymptotic discriminant}   of  $\Phi$   is the 
the   tail equivalence class,  $\cI(\Phi) = [\cI(\cU,x)]_{\infty}$,    for a   choice of the basepoint $x \in \fX$ and a choice of  an  adapted neighborhood basis $\cU$ at $x$. 
\end{defn}

   The explanation for the notation ``asymptotic discriminant'' will be made clear in Sections~\ref{subsec-stability} and \ref{subsec-ADGC}, where we give an interpretation of this class in terms of the group chain model for the action and the notion of the discriminant group for such actions, as  introduced in the works \cite{Dyer2015,DHL2016a,DHL2016c}.

Asymptotic discriminant defines an invariant of the return equivalence class of a pseudogroup action.
 \begin{prop}\label{prop-adaptedinv3}
 The asymptotic discriminant  $\cI(\Phi)$ of an equicontinuous minimal pseudogroup action $\Phi \colon  \cG \times \fX \to \fX$ on a Cantor space $\fX$ depends only on its return equivalence class.
\end{prop}
\proof
The proof is similar to that of Proposition~\ref{prop-retequivstable} so we give only a sketch of the argument.

Let    $\Phi^i \colon \cG_i \times \fX_i \to \fX_i$ be    equicontinuous  minimal pseudogroup Cantor actions,  for $i=1,2$, and assume that $\Phi^1$ and $\Phi^2$ are  return equivalent. Then are non-empty open sets $W_i \subset \fX_i$ and a homeomorphism $h \colon W_1 \to W_2$ such that the induced map $\Phi_{h} \colon \cG^*_{W_1} \to \cG^*_{W_2}$ is an isomorphism. 

Choose $x \in \fX_1$ and let $\cU = \{U_{\ell} \subset \fX_1  \mid \ell \geq 1\}$ be an  adapted neighborhood basis  at $x$ for $\Phi^1$. Let $y = \phi(x) \in \fX_2$ and $V_{\ell} = \phi(U_{\ell})$ for $\ell \geq 1$. Then $\cV = \{V_{\ell} \subset \fX_2 \mid \ell \geq 1\}$ is an  adapted neighborhood basis  at $y$ for $\Phi^2$. 

 The restricted homeomorphisms  $h_{\ell} = \phi | U_{\ell} \colon U_{\ell} \to V_{\ell}$ induce  isomorphisms
$h_{\ell}^* \colon \cH_{U_{\ell}} \to \cH_{V_{\ell}}$. By   Lemma~\ref{lem-closure-returnequiv} there are induced isomorphisms
$\overline{h_{\ell}} \colon \overline{\cH_{U_{\ell}}} \to \overline{\cH_{V_{\ell}}}$ of their closures in the uniform topology which coincides with the compact-open topology. 
Restrictions of these isomorphisms to isotropy groups gives the isomorphisms
$\overline{h_{\ell}} \colon \cI(U_1, U_{\ell},x) \to \cI(V_1, V_{\ell},x)$. Thus the maps $\{ \overline{h_{\ell}} \mid \ell \geq 1\}$ give an isomorphism of the asymptotic discriminant class associated to the adapted neighborhood basis  $\cU$ with the asymptotic class of the adapted neighborhood basis  $\cV$.  
\endproof

Let $\fM$ be an equicontinuous matchbox manifold, and choose a 
 regular covering  of $\fM$ as in Section~\ref{subsec-fs}, and let $\cGF$ be the associated holonomy pseudogroup as defined in Section~\ref{subsec-holonomy}, which yields an equicontinuous minimal pseudogroup action $\Phi \colon \cGF \times \fX \to \fX$ on the transversal Cantor space $\fX$. Let $x \in \fX$ and let $\cU = \{U_{\ell} \subset \fX  \mid \ell \geq 1\}$ be an  adapted neighborhood basis  at $x$ for $\Phi$. Define the asymptotic discriminant of $\fM$ to be  the asymptotic class $\cI(\fM) = \cI(\cU, x)$ associated to the adapted neighborhood basis $\cU$ at $x$. Then by   Theorem~\ref{thm-topinv} and Proposition~\ref{prop-adaptedinv3},      we have:

  \begin{thm}\label{thm-AD}
 The  asymptotic discriminant  $\cI(\fM)$ of $\fM$ is well-defined, and is an invariant of  the homeomorphism class of $\fM$.
  \end{thm}
  
The following result then follows from Definition~\ref{def-stable2} and Lemma~\ref{lem-asymptoticconstant}.
 \begin{cor}\label{cor-stable3}
 An equicontinuous matchbox manifold $\fM$ is stable if an only if the asymptotic discriminant invariant is asymptotically constant.
 \end{cor}
 
In the following sections, we show how to calculate the asymptotic discriminant for weak solenoids using the calculus of group chains, as developed in  \cite{Dyer2015,DHL2016a}.  The works \cite{Dyer2015,DHL2016a,DHL2016b,DHL2016c} used these methods to construct examples of various classes   of   weak solenoids which are stable and have non-trivial asymptotic discriminant. In fact,  Theorem~10.8 in  \cite{DHL2016c} shows that given any separable profinite group $\bK$ there is a weak solenoid with constant asymptotic discriminant isomorphic to $\bK$. This provides uncountable families of non-homeomorphic stable equicontinuous matchbox manifolds. 

Our interest in this work is on the equicontinuous matchbox manifolds which are not stable, hence are wild. For such spaces, the asymptotic discriminant is not a compact subgroup of a profinite group, but a tail equivalence class of maps between such groups. In the following, we introduce the methods used to calculate the asymptotic discriminant, then in Section~\ref{sec-wildex} give systematic construction of an uncountable collection of wild Cantor actions with distinct asymptotic discriminant groups.

\section{Weak solenoids}\label{sec-solenoids}

Sections~\ref{sec-ellis} and \ref{sec-adi} formulated the stable property and defined the asymptotic discriminant invariant 
in terms of holonomy pseudogroups, in analogy with the approach in  \cite{ALM2016}.  On the other hand, given an  equicontinuous minimal pseudogroup action, the works \cite{ClarkHurder2013,Dyer2015,DHL2016a}  associate to an equicontinuous minimal Cantor action    a group chain model for the action, and formulate many concepts and results analogous to those in the above sections in terms of the group chain model.   

We   discuss the group chain model associated to weak solenoids in this section, then   give an interpretation of the stable property and the asymptotic discriminant invariant  in terms of group chains in the following section. This is the basis for the construction of families of weak solenoids which are stable with non-trivial discriminant in the works \cite{DHL2016a,DHL2016b,DHL2016c}, and of families of weak solenoids which are wild  as described in Section~\ref{sec-wildex}.

We begin by recalling   the construction  of  \emph{weak  solenoids}, as first introduced by McCord \cite{McCord1965} and Schori \cite{Schori1966},   and       some of their properties as developed by Rogers and Tollefson \cite{Rogers1970,RT1971a,RT1971b,RT1972} and Fokkink and Oversteegen \cite{FO2002}.

  A  \emph{presentation} for a weak solenoid   is a collection 
\begin{equation}\label{eq-defweaksol}
\cP = \{ p_{\ell+1} \colon M_{\ell+1} \to M_{\ell} \mid \ell \geq 0\} \ ,
\end{equation}
    where each $M_{\ell}$ is a connected compact manifold   of dimension $n$, and each  \emph{bonding} map $p_{\ell +1}$  is a proper covering map of finite index.
    Associated to a presentation $\cP$ is the weak solenoid denoted by $\cS_{\cP}$ which is defined as the inverse limit,  
\begin{equation}\label{eq-presentationinvlim}
\cS_{\cP} \equiv \lim_{\longleftarrow} ~ \{ p_{\ell +1} \colon M_{\ell +1} \to M_{\ell}\} ~ \subset \prod_{\ell \geq 0} ~ M_{\ell} ~ .
\end{equation}
 By definition, for a sequence $\{x_{\ell} \in M_{\ell} \mid \ell \geq 0\}$, we have 
\begin{equation}\label{eq-presentationinvlim2}
x = (x_0, x_1, \ldots ) \in \cS_{\cP}   ~ \Longleftrightarrow  ~ p_{\ell}(x_{\ell}) =  x_{\ell-1} ~ {\rm for ~ all} ~ \ell \geq 1 ~. 
\end{equation}
The set $\cS_{\cP}$ is given  the relative  topology, induced from the product topology, so that $\cS_{\cP}$ is itself compact and connected.
  McCord showed in \cite{McCord1965} that  the space $\cS_{\cP}$ has a local product structure, and moreover we have:

\begin{prop}\label{prop-solenoidsMM}
Let  $\cP$ be a presentation with base space $M_0$ and  $\cS_{\cP}$ the associated weak solenoid. Then   $\cS_{\cP}$ is  an equicontinuous matchbox manifold of dimension $n$ with foliation   $\F_{\cP}$  by path-connected components. 
\end{prop}

Associated to a presentation $\cP$ of compact manifolds is a sequence of proper surjective maps 
\begin{equation}\label{eq-coverings}
q_{\ell} = p_{1} \circ \cdots \circ p_{\ell -1} \circ p_{\ell} \colon M_{\ell} \to M_0 ~ .
\end{equation}
For each $\ell > 1$, projection onto the $\ell$-th factor in the product $\ds \prod_{\ell \geq 0} ~ M_{\ell}$ in \eqref{eq-presentationinvlim} yields a 
  fibration map denoted by $\Pi_{\ell} \colon \cS_{\cP}  \to M_{\ell}$, for which 
 $\Pi_0 = \Pi_{\ell} \circ q_{\ell} \colon \cS_{\cP} \to M_0$. 
A choice of a basepoint $x_0 \in M_0$ fixes a fiber  $\fX_0 = \Pi_0^{-1}(x_0)$, which is a Cantor set by the assumption  that the fibers of each map $p_{\ell}$ have cardinality at least $2$.  The choice of $x_0$ will remain fixed throughout the following. We also then have a fixed base group   $G_0 =  \pi_1(M_{0}, x_{0})$.

A choice  $x \in \fX_0$   defines basepoints  $x_{\ell} = \Pi_{\ell}(x) \in M_{\ell}$ for $\ell \geq 1$.
For each $\ell \geq 1$, let  
\begin{equation}\label{eq-imahes}
G^x_{\ell} = {\rm image}\left\{  (q_{\ell} )_{\#} \colon \pi_1(M_{\ell}, x_{\ell}) \longrightarrow   G_0\right\}  
\end{equation}
  denote  the image of the induced map $(q_{\ell} )_{\#} $ on fundamental groups. Associated to the presentation $\cP$ and basepoint $x \in \fX_0$ we thus obtain a descending chain of subgroups of finite index
  \begin{equation}\label{eq-descendingchain}
\cG^x \equiv  \{ G^x_{\ell} \}_{\ell \geq 0} = \left\{G_0 \supset G^x_{1} \supset G^x_{2} \supset \cdots \supset G^x_{\ell} \supset \cdots \right\} \ .
\end{equation}

Note that given another choice of basepoint $y \in \fX_0$ with corresponding images $y_{\ell} =   \Pi_{\ell}(y) \in M_{\ell}$, for each $\ell \geq 1$, there exists $g^y_{\ell} \in G_0$ such that the image group $G^y_{\ell} = (g^y_{\ell}) G^x_{\ell} (g^y_{\ell})^{-1}$. The resulting group chain $\cG^y = \{ G^y_{\ell} \}_{\ell \geq 0}$ is said to be \emph{conjugate} to $\cG^x$.

Each quotient  $X_{\ell}^x = G_0/G_{\ell}^x$ is   a   finite set equipped with a left $G_0$-action, and there are surjections $X_{\ell +1}^x \to X_{\ell}^x$ which commute with the action of $G_0$.  The inverse limit 
\begin{equation}\label{eq-Galoisfiber}
X_{\infty}^x = \lim_{\longleftarrow} ~ \{ p_{\ell +1} \colon X_{\ell +1}^x \to X_{\ell}^x\} ~ \subset \prod_{\ell \geq 0} ~ X_{\ell}^x  
\end{equation}
is given the relative topology, induced from   the product (Tychonoff) topology on the   space $\ds \prod_{\ell \geq 0} ~ X_{\ell}^x$, so that $\ds X_{\infty}^x$ is compact. 
As each   set $X_{\ell}^x$ is finite with    cardinality at least $2$ for $\ell \geq 1$, $X_{\infty}^x$  is   a totally disconnected perfect set, so is   a Cantor space. 

A sequence $(g_{\ell})  \subset G_0$ such that $g_{\ell} G_{\ell}^x = g_{\ell+1} G_{\ell}^x$ for all $\ell \geq 0$ determines a point $(g_{\ell} G_{\ell}^x) \in X_{\infty}^x$. 
Let $e \in G_0$ denote the identity element, then the sequence $e_{x} = (e G_{\ell}^x)$ is the basepoint of $X_{\infty}^x$.

The action  $\Phi_x \colon G_0 \times X_{\infty}^x \to X_{\infty}^x$ is given by coordinate-wise multiplication, $\Phi_0(g)(g_{\ell} G_{\ell}^x) = (g g_{\ell} G_{\ell}^x)$.
 We then have the standard observation:

 \begin{lemma}\label{lem-denseaction}
 $\Phi_x \colon G_0 \times X_{\infty}^x \to X_{\infty}^x$ defines an equicontinuous     Cantor minimal system  $(X_{\infty}^x , G_0 , \Phi_x)$.  
\end{lemma}  
 
The choice of the basepoint $x \in \cS_{\cP}$ defines   basepoints $x_{\ell} \in M_{\ell}$ for all $\ell \geq 1$, which gives   an identification of $X_{\ell}^x$ with the fiber of the covering map $M_{\ell} \to M_0$. In the inverse limit, we thus obtain a homeomorphism $\tau_x \colon X_{\infty}^x \to  \fX_0 = \Pi_0^{-1}(x_0)$ such that $\tau_x(e_x) = x$, 
which can be viewed as ``coordinates'' on $\fX_0$.

 The left action of $G_0$ on $X_{\infty}^x$ is conjugated to an action of $G_0$ on $\fX_0$,  called the \emph{monodromy action} at $x_0$ for the fibration $\Pi_0 \colon \cS_{\cP} \to M_0$, where the   action is  defined by the holonomy transport along the leaves of the foliation $\F_{\cP}$ on $\cS_{\cP}$.

The map  $\tau_x \colon X_{\infty}^x \to  \fX_0$ is used to give   a ``standard form'' for the solenoid $\cS_{\cP}$.
Let $\wtM_0$ denote the universal covering of the compact manifold $M_0$ and let $(X_{\infty}^x , G_0 , \Phi_x)$ be the minimal Cantor system associated to the presentation $\cP$ and the choice of a basepoint $x \in \fX_0$. Associated to the left action $\Phi_x$ of $G_0$ on $X_{\infty}^x$ is a suspension space 
\begin{equation}\label{eq-suspensionfols}
\fM = \wtM_0 \times X_{\infty}^x / (z \cdot g^{-1}, x) \sim (z , \Phi_x(g)(x)) \quad {\rm for }~ z \in \wtM_0 , ~ g \in G_0 ~,
\end{equation}
which is a minimal matchbox manifold.  Moreover, $\fM$ has an inverse limit  presentation,  where all of the bonding maps between the coverings $M_{\ell} \to M_0$ are derived from the universal covering map $\wtpi \colon \wtM_0 \to M_0$, so are in ``standard form''. The following result uses the path lifting property for coverings, to show that for any presentation $\cP$, we have:
\begin{thm} \cite{ClarkHurder2013}\label{thm-weaksuspensions}
Let   $\cS_{\cP}$ be a weak solenoid, with base space $M_0$ which is a compact manifold of dimension $n \geq 1$. Then the suspension of the map $\tau_x$ yields a foliated homeomorphism $\tau_x^* \colon \fM \to \cS_{\cP}$.
\end{thm}
 
  \begin{cor}\label{cor-weaksuspensions}
The homeomorphism type of a weak solenoid    $\cS_{\cP}$ is completely determined by the base manifold $M_0$ and the associated minimal Cantor system   $(X_{\infty}^x , G_0 , \Phi_x)$.
\end{cor}

 We conclude this discussion of  this section, by introducing  the following important notion:
 \begin{defn}\label{def-kernel}
The \emph{kernel} of a group chain  $\cG = \{G_{\ell}\}_{\ell \geq 0}$ is the subgroup  $\ds K(\cG) = \bigcap_{\ell \geq 0} ~ G_{\ell}$.
\end{defn}
  
For a weak solenoid $\cS_{\cP}$ with choice   of a basepoint $x_0 \in M_0$ and  fiber  $\fX_0 = \Pi_0^{-1}(x_0)$, 
  the kernel subgroup $K(\cG^x) \subset G_0$  may depend on the choice of the basepoint $x \in \fX_0$. 
  The dependence   of $K(\cG^x)$ on $x$ is a natural aspect  of the dynamics of the foliation $\F_{\cP}$ on $\cS_{\cP}$, when 
  $K(\cG^x)$ is interpreted in terms of the topology of the leaves of $\F_{\cP}$.
  
The map $\tau_x^* \colon \fM \to \cS_{\cP}$ of Theorem~\ref{thm-weaksuspensions} sends the 
  quotient space  $\wtM/K(\cG^x)$    to the leaf $L_x \subset \cS_{\cP}$ through $x \in \fX_0$ in $\cS_{\cP}$ and so $K(\cG^x)$  is naturally identified with the fundamental group $\pi_1(L_x , x)$.
 The holonomy homomorphism $h_x \colon \pi_1(L_x , x) \to Homeo(\fX_0 , x)$    of the leaf $L_x$ in the suspension foliation    $\F_{\cP}$ is   conjugate to the left action,  $\Phi_0 \colon K(\cG^x) \to Homeo(X_{\infty}^x , e_x)$.

 From the point of view of foliation theory, the leaves of $\F_{\cP}$ with holonomy are a ``small'' set''.  There   always exists leaves without holonomy, while there may exist leaves with holonomy, and so the fundamental groups of the leaves may vary accordingly.  This aspect of the foliation dynamics of weak solenoids is discussed further in \cite[Section~4.2]{DHL2016c}.

\section{Algebraic models from group chains}\label{sec-stablechains}
In this section, we develop models derived from a group chain $\cG$ for the spaces and actions introduced in Section~\ref{sec-ellis}. In particular, we construct  
the Ellis group for the equicontinuous action  on the Cantor space  $X_{\infty}^x$  associated to a group chain in \eqref{eq-Galoisfiber}, and then interpret in Section~\ref{subsec-stability} the definition of stable actions given in Section~\ref{subsec-defSQA}  in terms of the group chain model of the action. 
In Section~\ref{subsec-ADGC}, we then give   the group chain interpretation   of the asymptotic discriminant invariant in Definition~\ref{def-asympdisc}.

\subsection{Core chains}\label{subsec-cores}
Let    $\cG  = \{G_{\ell}\}_{\ell \geq 0}$ be a group chain in $G_0$. Let $X_{\infty}$ be the inverse limit space defined as in \eqref{eq-Galoisfiber} by the finite quotients  $X_{\ell} = G_0/G_{\ell}$ with the transitive left $G_0$-action. This defines an equicontinuous minimal Cantor action denoted by 
  $(X_{\infty} , G_0 , \Phi_0)$.   Let $\cH = \Phi_0(G_0) \subset Homeo(X_{\infty})$ be the subgroup defined by this action,  and let $\ocH \subset Homeo(X_{\infty})$ be its closure.

  For each $\ell \geq 0$,   the group $G_0$ acts transitively on the left on each quotient space $X_{\ell}$, 
  where $G_{\ell}$ is the isotropy group for the basepoint $x_{\ell} = e G_{\ell} \in X_{\ell}$. 
   Then for $g \in G_0$ the isotropy of the point $g G_{\ell} \in X_{\ell}$ is the conjugate subgroup $g G_{\ell} g^{-1} \subset G_0$. 
   The kernel of the homomorphism $\Phi_{\ell} \colon G_0 \to Homeo(X_{\ell})$ is the intersection of the isotropy groups for the points of $X_{\ell}$, 
   \begin{equation}\label{eq-core}
C_{\ell} \equiv   {\rm core}_{G_0}   \, G_{\ell} \equiv  \bigcap_{g \in {G_0}} \ gG_{\ell} g^{-1}  ~ \subseteq ~ G_{\ell} \ .
\end{equation}
The subgroup $C_{\ell}$, called the     \emph{core} of $G_{\ell}$ in $G_0$, is the maximal normal subgroup of $G_{\ell}$.   Note that 
\begin{equation}\label{eq-normkernel}
\bigcap_{\ell \geq 0} \ C_{\ell} = \bigcap_{\ell \geq 0} \  \bigcap_{g \in {G_0}} \  gG_{\ell} g^{-1}  
=  \bigcap_{g \in {G_0}} \   \bigcap_{\ell \geq 0} \  gG_{\ell} g^{-1} 
=   \bigcap_{g \in {G_0}} \  g K(\cG) g^{-1}  \equiv N(\cG)   \ ,
\end{equation}
so that  $N(\cG) \subset K(\cG)$ is the largest normal subgroup of $G_0$ contained in the kernel $K(\cG)$.

Each  quotient $G_0/C_{\ell}$ is a finite group, and  the collection  $\cC = \{C_{\ell} \}_{\ell \geq 0}$ forms a descending chain of normal subgroups of $G_0$. The inclusions of coset spaces define bonding homomorphisms $\delta^{\ell+1}_{\ell} $ for the inverse sequence of quotient groups  $G_0/C_{\ell}$, and the inverse limit space 
\begin{eqnarray}
  C_{\infty}    & = &  \lim_{\longleftarrow} \, \left\{\delta^{\ell+1}_{\ell} \colon   G_0/C_{\ell+1} \to G_0/C_{\ell}   \right\}  \label{cinfty-define}  \\
  & = &  \{(eG_0, g_1 C_1, \ldots) \mid g_{\ell} C_{\ell} = g_{\ell+1} C_{\ell}\} ~ \subset \prod_{\ell \geq 0} \ G_0/C_{\ell}  \label{cinfty-coords} 
\end{eqnarray} 
is a profinite group. Let  $\whiota \colon G_0 \to C_{\infty}$ be the homomorphism defined by $\whiota(g) = (gC_{\ell})$ for $g \in G_0$. 
Note that $N(\cG) = ker(\whiota)$, so that $\whiota$ is injective if and only if the kernel $K(\cG)$ has trivial core. 

Let $
\whPhi_0$ be  the  induced left action of $G_0$ on $C_{\infty}$, where 
$\whPhi_0(g)   (g_{\ell} C_{\ell}) = (gg_{\ell} C_{\ell})$ for $g \in G_0$.
  Let  $(C_{\infty} , G_0, \whPhi_0)$ denote the minimal Cantor system defined by this action. As $C_{\infty}$ is a group, the isotropy subgroup of   the identity $\whe_{\infty} = (eC_{\ell}) \in C_{\infty}$ is the   subgroup $N(\cG) \subset G_0$ which is the kernel of the homomorphism $\Phi_0 \colon G_0 \to Homeo(C_{\infty})$.

Introduce the descending chain of clopen neighborhoods of the identity $\whe_{\infty}$ in $C_{\infty}$:
\begin{eqnarray}  
C_{n,\infty}   & = &  \lim_{\longleftarrow} \, \left\{\delta^{\ell+1}_{\ell} \colon   G_n/C_{\ell +1}  \to G_n/C_{\ell}  \mid \ell \geq n \right\}   \label{Cinfty-definenn} \\
  & = & \{ (g_{\ell} C_{\ell}) \in C_{\infty} \mid   g_\ell \in G_n   \} \ . \nonumber \ 
 \end{eqnarray}
 Note that $C_{n,\infty}$ can also be identified with the closure in $C_{\infty}$ of the image of $\whiota(G_n)$.

The collection $\{C_{n,\infty}  \mid n \geq 1\}$  of clopen subsets of $C_{\infty}$ defines a neighborhood system of $\whe_{\infty}$.

Observe that for each $\ell \geq  0$, the    quotient group $D_{\ell}  = G_{\ell}/C_{\ell} \subset G_0/C_{\ell}$ and so the inverse limit   
\begin{equation}\label{eq-discriminantdef}
\cD  = \lim_{\longleftarrow}\, \left \{\delta^{\ell+1}_{\ell} \colon  D_{\ell+1}  \to D_{\ell}  \right\}
\end{equation}
is a closed subgroup of $C_\infty$. The group $\cD$ is called the \emph{discriminant group} of the action $(V_0,G_0,\Phi_0)$.

The   relationship between $C_\infty$ and the Ellis group of   $(X_{\infty} , G_0 , \Phi_0)$ is given by the following result.

\begin{thm}[Theorem~4.4, \cite{DHL2016a}]\label{thm-quotientspace}
Let   $(X_{\infty} , G_0 , \Phi_0)$ be the   equicontinuous minimal Cantor action associated to the group chain $\cG  = \{G_{\ell}\}_{\ell \geq 0}$. 
Then the map $\whiota \colon G_0 \to C_{\infty}$ induces an isomorphism
 of topological groups $\Upsilon_0 \colon \overline{\Phi_0(G_0)} \to C_{\infty}$ such that the restriction gives an isomorphism $\Upsilon_{\whe} \colon \overline{\Phi_0(G_0)}_{\whe} \cong \cD$.
\end{thm}
 \proof
 We give the essential idea of the proof. For $\ell \geq 1$, we have $N(\cG) \subset C_{\ell}$, so the image $\Phi_0(G_0) \cong G_0/N(\cG)$ maps onto the finite quotient group $G_0/C_{\ell}$. Moreover, the  subgroup    $G_{\ell}/C_{\ell}$ is  the isotropy group of the action of $G_0/C_{\ell}$ on $X_{\ell}$.
 Passing to   inverse limits,  we obtain the isomorphisms $\Upsilon_0 \colon \overline{\Phi_0(G_0)} \to C_{\infty}$ and $\Upsilon_{\whe} \colon \overline{\Phi_0(G_0)}_{\whe} \to \cD$. 
  \endproof

The following is a key observation about the discriminant subgroup defined by \eqref{eq-discriminantdef}.

\begin{prop}[Proposition~5.3, \cite{DHL2016a}]\label{prop-discrisnotnormal}
Let   $(X_{\infty} , G_0 , \Phi_0)$ be the   equicontinuous minimal Cantor action associated to the group chain $\cG  = \{G_{\ell}\}_{\ell \geq 0}$. 
Then  
\begin{equation}\label{eq-rationalcore}
 {\rm core}_{G_0} \  \cD  ~ = ~    \bigcap_{g \in G_0} \ g \ \cD \ g^{-1}  ~ = ~ \{\whe\} \ .
\end{equation}
Moreover, the   maximal normal subgroup in $C_{\infty}$ of $\cD$ is trivial.
\end{prop}

\subsection{Stable group chains}\label{subsec-stability}
We next consider    the stability properties of the discriminant  group  for the      minimal Cantor system  $(X_{\infty} , G_0 , \Phi_0)$   associated to the group chain $\cG  = \{G_{\ell}\}_{\ell \geq 0}$.  We use the notations of Section~\ref{subsec-defSQA}. For $g \in G_0$ and $(g_{\ell} G_{\ell}) \in X_{\infty}$ let $g \cdot (g_{\ell}G_{\ell}) = \Phi_0(g)(g_{\ell}G_{\ell})$. 

For $n \geq 1$, define the clopen   neighborhood   of $\whx_0 = (e G_{\ell}) \in X_{\infty}$, 
\begin{equation}  \label{Ginfty-definen}
U_n       =   \{ (g_{\ell} G_{\ell}) \in X_{\infty} \mid   g_\ell \in G_n   \} \subset X_{\infty}  \ .
\end{equation}

Note that for each $n \geq 1$,  the $G_0$-translates of the set $U_n$ form a partition of  $X_{\infty}$ so the collection $\ds \{g \cdot U_n \mid n \geq 1 ~ , ~ g \in G_0\}$ forms a basis of clopen sets for the topology of $X_{\infty}$. In particular, 
 the restriction of the $G_0$-action on $X_{\infty}$ to the set $U_n$ is given by 
the  action of $G_n$ on $U_n$.

Let $\cH_{U_n} = \Phi_0(G_n) \subset Homeo(U_n)$, then the Ellis group for the restricted action to $U_n$  is the closure 
  $\overline{\cH_{U_n}} = \overline{\Phi_0(G_n)}  \subset Homeo(U_n)$. 
   For $m > n  \geq 0$,  note that    $\cH_{U_n,U_m}  = \Phi_0(G_m) \subset \cH_{U_n}$, 
  and so $\overline{\cH_{U_n,U_m}} = \overline{\Phi_0(G_m)} \subset \overline{\cH_{U_n}}$. 
  We  next express   the restriction map $\orho_{U_n,U_m} \colon \overline{\cH_{U_n,U_m}} \to \overline{\cH_{U_m}}$   and its kernel in terms of group chains and their properties.

We first derive  an inverse limit   presentation for $\overline{\cH_{U_n}}$.
 Introduce the truncated group chain associated to the restricted action of $G_n$ on $U_n$, 
 \begin{equation}\label{eq-truncatedchain}
\cG_n = \left\{ G_{\ell} \right\}_{\ell \geq n} = \left\{ G_n \supset G_{n+1} \supset G_{n+2} \supset     \cdots \right\} \ .
\end{equation}
 For each $\ell \geq n \geq 0$, we   have the core subgroup  for the truncated chain $\cG_n$ defined by 
 \begin{equation}\label{eq-Mn}
C_{n,\ell} \equiv   {\rm core}_{G_n}   \, G_{\ell}  \equiv  \bigcap_{g \in {G_n}} gG_{\ell} g^{-1}   \subset G_n \ .
\end{equation}
Observe that $C_{n,\ell}$ is the kernel of the action of $G_n$ on the quotient set $G_n/G_{\ell}$.
Moreover,  $C_{0,\ell}= C_{\ell}$, and   for all $\ell > m \geq n \geq 0$, we have $\ds C_{n,\ell}  \subset C_{m,\ell} \subset G_{\ell} \subset G_m \subset G_n$. Define the profinite group 
    \begin{eqnarray}  
\fC_{n,m}   & \cong &  \lim_{\longleftarrow} \, \left\{\delta^{\ell+1}_{\ell} \colon   G_{m}/C_{n,\ell} \to G_m/C_{n,\ell + 1}  \mid \ell \geq m  \right\}   \\
   & = & \{ (g_{\ell} C_{n,\ell})   \mid  \ell \geq m   \ , \  g_{m} \in G_{m} \ , \  g_{\ell +1} C_{n,\ell} = g_{\ell} C_{n,\ell}  \}     \label{Dinfty-definen}  \ .  
   \end{eqnarray}
Then $\fC_{0,0} = C_{\infty}$,   $\fC_{n,n}$ is the Ellis   group   for the  truncated group chain $\cG_n$ and 
by definition  we have that $\fC_{n,m} \subset \fC_{n,n}$.  The group $\fC_{n,m} $ is a clopen neighborhood of the identity in $\fC_{n,n}$. 
For example, $\fC_{0,m} =C_{m,\infty}$ in \eqref{Cinfty-definenn}.

For $n \geq 0$, the  homomorphism $\iota_n \colon G_n \to \fC_{n,n}$ has dense image, hence induces    an isomorphism $\Upsilon_n \colon \overline{\cH_{U_n}} \to \fC_{n,n}$, which follows by    analogous reasons as that of Theorem~\ref{thm-quotientspace}. Thus $\overline{\cH_{U_n}} = \fC_{n,n}$ is the Ellis group of the action of $G_n$ on $U_n$.

 Next, for $m \geq n \geq 0$, the subgroup $\iota_n(G_m) \subset \fC_{n,m}$ is dense, hence induces    an isomorphism $\Upsilon_{n,m} \colon  \overline{\cH_{U_n,U_m}} \to \fC_{n,m} $. 
 That is, $\fC_{n,m} \subset \fC_{n,n}$ is a group chain model for     $\overline{\cH_{U_n,U_m}} \subset \overline{\cH_{U_n}}$.

   For each $\ell \geq k \geq m \geq n$, the inclusions $C_{n,\ell} \subset C_{m,\ell}$ induce   group surjections  denoted by 
\begin{align}\label{eq-quotiensmn} 
\phi_{k,n,m}^{\ell} \colon   G_{k}/C_{n,\ell}  \longrightarrow  G_{k}/C_{m,\ell} \ .
\end{align}
 For   $k \geq m \geq n \geq 0$, standard methods as in \cite{RZ2010} show that this yields   surjective homomorphisms  of profinite groups 
   $\ds \whphi_{k,n,m} \colon   \fC_{n,k}  \to \fC_{m,k}$ which commute with the left action of $G_0$.
In particular, for $k=m$,    we obtain a surjective homomorphism  
  \begin{align}\label{eq-chomeomorph} \whphi_{m,n,m} \colon   \fC_{n,m}  \to \fC_{m,m}\end{align} 
  which is a chain representation of the restriction map  $\orho_{U_n,U_m}$. We   express its kernel   in terms of the discriminant groups of the truncated group chains $\cG_n$.

For $n \geq 0$,   the  discriminant group  associated to the group chain   $\cG_n$ is given by,  for  $m \geq n$, 
\begin{eqnarray} 
\cD_n & = &    \lim_{\longleftarrow}\, \left \{\delta^{\ell+1}_{\ell} \colon  G_{\ell+1}/C_{n,\ell+1}  \to G_{\ell}/C_{n,\ell} \mid \ell \geq n\right\}    \subset \fC_{n,n} \label{eq-discquotients1} \\
  & \cong &    \lim_{\longleftarrow}\, \left \{\delta^{\ell+1}_{\ell} \colon  G_{\ell+1}/C_{n,\ell+1}  \to G_{\ell}/C_{n,\ell} \mid \ell \geq m \right\}  \subset \fC_{n,m} \ . \label{eq-discquotients2}
\end{eqnarray}
Let $\cD_{n,m} \subset \fC_{n,m}$ denote the image of $\cD_n$ under the map \eqref{eq-discquotients2}.
 It then follows from  \eqref{eq-quotiensmn}    that for $m > n \geq 1$,  there are    surjective homomorphisms, 
\begin{equation}\label{eq-discmapsnm2}
   \cD   ~ \stackrel{~ \psi_{0,n} ~ }{\longrightarrow} ~  \cD_n \cong \cD_{n,m} ~ \stackrel{~ \psi_{n,m} ~}{\longrightarrow} ~  \cD_m \ ,
\end{equation}
 where the map $\psi_{n,m}$ is the restriction of the surjection  $\ds \whphi_{m,n,m} \colon   \fC_{n,m}  \to \fC_{m,m}$ in \eqref{eq-chomeomorph}.

  We can now formulate the notion of a stable group chain, as introduced in the work \cite{DHL2016c}.
\begin{defn}\label{def-stableGC}
A group chain $\cG = \{G_{\ell}\}_{\ell \geq 0}$ is said to be \emph{stable} if there exists $n_0 \geq 0$ such that for all $m \geq n \geq n_0$, the map  $\ds \psi_{n,m} \colon \cD_n  \to    \cD_m$ defined in \eqref{eq-discmapsnm2} is an isomorphism. Otherwise, the group chain is said to be \emph{wild}. 
\end{defn}

Here is one of our main results, relating the results in \cite{DHL2016c} with the notion of stability in Section \ref{subsec-defSQA}.
\begin{prop}\label{prop-stableisstable}
Let $\cG = \{G_{\ell}\}_{\ell \geq 0}$ be a group chain. Then the     equicontinuous minimal Cantor action   $(X_{\infty}, G_0 , \Phi_0)$ associated to $\cG$  is stable in the sense of Definition~\ref{def-stable2} if and only if the group chain $\cG$ is stable in the sense of Definition~\ref{def-stableGC}.
\end{prop} 
 \proof
Recall that  for $U_m \subset X_{\infty}$ defined by \eqref{Ginfty-definen} we have  $G_m = \{g \in G_0 \mid \Phi_0(g)( U_m) = U_m\}$.
The   group $\fC_{m,m}$ is the inverse limit group formed from the quotients $G_m/C_{m,\ell}$ where each $C_{m,\ell} \subset G_m$ for $\ell \geq m$ is normal in $G_m$, and 
acts transitively on   $U_m$ with isotropy group $\cD_m$. 
The group  $ \fC_{n,m}$ is the inverse limit group formed from the quotients $G_m/C_{n,\ell}$ where each $C_{n,\ell} \subset G_\ell \subset G_m$   for $\ell \geq m$ is normal in $G_n$, hence is also normal in $G_m$. The group  $ \fC_{n,m}$ 
 acts transitively on   $U_m$ with isotropy group   $\cD_{n,m}$. 
Thus, we have the commutative diagram of fibrations,
    \begin{align}\label{eq-diagramkernel} 
    \xymatrixcolsep{4pc}
    \xymatrixrowsep{1.6pc}
   \xymatrix{
   \cD_{n,m}  \ar[d]   \ar[r]^{\psi_{n,m}} & \cD_m \ar[d]  \\
   \fC_{n,m} \ar[d] \ar[r]^{\whphi_{m,n,m}}  & \fC_{m,m} \ar[d]   \\
 U_m  \ar[r]^{=} & U_m   
   }  
   \end{align}

   Assume first that the group chain $\cG$ is stable  in the sense of Definition~\ref{def-stableGC}, and let $n_0 \geq 0$ be such  that $\ds \psi_{n,m} \colon \cD_{n,m}  \to    \cD_m$  is an isomorphism  for all $m \geq n \geq n_0$. Then    the center  map $\whphi_{m,n,m}$ in \eqref{eq-diagramkernel} is an isomorphism.    
   Thus, for $\whg \in \fC_{n,m}$ if  $\whphi_{m,n,m}(\whg) \in \fC_{m,m}$ acts trivially on $U_m$ then $\whg$ also acts trivially on $U_m$. That is, the   restriction map  $\orho_{U_n,U_m}$ is injective. As this holds for all $m > n$, by Proposition~\ref{prop-retequivstable}  this implies that the action $\Phi_0$ is stable  in the sense of Definition~\ref{def-stable2}.

  Now suppose that   the action $\Phi_0 \colon G_0 \times X_{\infty} \to X_{\infty}$ is stable  in the sense of Definition~\ref{def-stable2}.  The collection of clopen sets $\{U_{\ell} \mid \ell \geq 0\}$ is a neighborhood basis around $x_0 \in X_{\infty}$, so there exists   $n_0 > 0$ be such that   the restriction map $\orho_{U_{n},U_m}$  {on closures}   in \eqref{eq-restriction} is an isomorphism for all $m \geq n \geq n_0$.
We claim    that $\ds \psi_{n,m} \colon \cD_n \cong  \cD_{n,m}  \to    \cD_m$  is an isomorphism.

Let $\whg \in \overline{\cH_{U_n,U_m}}$, and suppose that   $\orho_{U_n,U_m}(\whg)$ is the identity. 
Then we can write $\whg = (g_{\ell} C_{n,\ell}) \in \fC_{n,m}$ where $g_{\ell} \in G_{m}$ for $\ell \geq m$, 
and  the element $\whphi_{m,n,m}(g_{\ell} C_{n,\ell}) = (g_{\ell} C_{m,\ell}) \in \fC_{m,m}$ acts as the identity by multiplication on the left on $U_m \cong \fC_{m,m}/\cD_m$.

That is, for all $(h_{\ell} C_{m,\ell}) \in \fC_{m,m}$,   we have the identity of cosets, 
$$(g_{\ell} C_{m,\ell}) \cdot (h_{\ell} C_{m,\ell}) \ \cD_m = (h_{\ell} C_{m,\ell}) \ \cD_m \ ,$$
which implies that 
  $(g_{\ell} C_{m,\ell})     \in (h_{\ell} C_{m,\ell}) \cD_m (h_{\ell} C_{m,\ell})^{-1}$ for all 
$(h_{\ell} C_{m,\ell}) \in \fC_{m,m}$.
Then  
\begin{equation}\label{eq-rationalcorek}
 {\rm core}_{G_m} \cD_m  =  \bigcap_{h \in G_m} h \ \cD_m \ h^{-1}  = \{\whe_m\} ~ , ~ \whe_m = (e    C_{m,\ell})  \ , 
\end{equation}
 implies that  $(g_{\ell} C_{m,\ell}) =   \whe_m$, or $g_{\ell} \in C_{m,\ell}$ for all $\ell \geq m$.
   Thus we may suppose     that   $\whg = (g_{\ell} C_{n,\ell})$ with $g_{\ell} \in C_{m,\ell}$ for all $\ell \geq m$.
   
By the assumption that the action is stable, we have that $\whg$ acts trivially on $U_n$. In particular, for all  $(h_{\ell} C_{n,\ell}) \in \fC_{n,m}$,   we have the identity of cosets, 
\begin{equation}\label{eq-identitynCn}
(g_{\ell} C_{n,\ell}) \cdot (h_{\ell} C_{n,\ell}) \ \cD_{n,m} = (h_{\ell} C_{n,\ell}) \ \cD_{n,m} \ ,
\end{equation}

and by a similar argument $(g_\ell C_{n,\ell}) = \whe_n $.

This shows that if $\whg \in \fC_{n,m}$ is mapped onto $\whe \in \fC_{m,m}$, then $\whg = \whe_n$, and the kernel of the map $\whphi_{m,n,m}$ is trivial. Since $\psi_{n,m}$ is the restriction of $\whphi_{m,n,m}$, it also has a trivial kernel. This shows that the group chain $\cG$ is stable in the sense of Definition~\ref{def-stableGC}.
 \endproof

\subsection{Asymptotic discriminant for group chains}\label{subsec-ADGC}

We   consider   the asymptotic  discriminant  invariant for    for a      minimal Cantor system  $(X_{\infty} , G_0 , \Phi_0)$   associated to the group chain $\cG  = \{G_{\ell}\}_{\ell \geq 0}$.  We use the notations of Section~\ref{subsec-defSQA} and Section~\ref{sec-adi}.

Recall the definition of the clopen sets $U_{n}$ in \eqref{Ginfty-definen}. It was observed that the restriction of the $G_0$-action on $X_{\infty}$ to the set $U_n$ is given by  the  action of $G_n$ on $U_n$, so that the collection of clopen sets   $\cU = \{U_{\ell} \mid \ell \geq 0\}$ is an adapted neighborhood basis at $x_0 \in X_{\infty}$ for the action $\Phi_0$.

The discriminant for the action of $G_n$ on $U_n$ is given by $\cD_n$ as in \eqref{eq-discquotients1}, which 
either by  Theorem~\ref{thm-quotientspace} or direct calculation, is identified with the isotropy group at $x_0$ of the action of $\fC_{n,n}$ on $U_n$.
Thus, in the notation of the definition   \eqref{eq-defisotropy}, we have the identification $\cD_n = \cI(U_1, U_{n},x_0)$. 
Moreover, the surjective maps  in Corollary~\ref{cor-restrictionsrho} are identified with the surjective maps   in  \eqref{eq-discmapsnm2}, so for $m > n \geq 1$ there is a commutative diagram
\begin{align} 
    \xymatrixcolsep{.1pc}
  \xymatrixrowsep{1.6pc}
 \xymatrix{
\sigma_{U_{n},U_{m}}   & \colon &   \cI(U_1, U_{n},x) \ar[d]_{\cong} & \longrightarrow &   \cI(U_1, U_{m},x)  \ar[d]^{\cong}     \\
 \psi_{n,m}  &   \colon &   \cD_n \cong \cD_{n,m} & \longrightarrow &    \cD_m  \cong \cD_{m,m}      
} 
\end{align}

 Then in analogy with  Definition~\ref{def-asympdisc}, we introduce the following invariant for $\cG$.
 
  \begin{defn}\label{def-asympdisc2}
 Let $\cG  = \{G_{\ell}\}_{\ell \geq 0}$ be a group chain. Then the \emph{asymptotic discriminant} for $\cG$ is the tail equivalence class $[\cD(\cG)]_{\infty}$ of   
 the sequence of surjective group homomorphisms 
 \begin{equation}\label{eq-surjectiveDn}
\cD(\cG) = \{\psi_{n,n+1} \colon \cD_n \to \cD_{n+1} \mid n \geq 1\}
\end{equation}
 defined by the discriminant groups for the restricted actions of $G_n$ on the clopen sets $U_n \subset \fX_{\infty}$.
  \end{defn}
The definition of $\cD(\cG)$ above involves a minor  abuse of notation, or at least a change in usage for the notation. In the works  \cite{Dyer2015,DHL2016a,DHL2016b} the discriminant group of a group chain was defined to be $\cD = \cD_1$ while   the work  \cite{DHL2016c} studied the   discriminant groups $\cD_n$ as a function of $n$, and was concerned mainly with the case when the group chain is stable, so that $\cD_n \cong \cD_m$ for some $n_0 \geq 1$ and $m > n \geq n_0$. 
By Lemma~\ref{lem-asymptoticconstant}, in the stable case the sequence of surjective group homomorphisms   $\cD(\cG)$ as defined above is asymptotically constant, so the two usages are thus compatible. For the case when $\cD(\cG)$ is unstable, or wild, then the notion of the discriminant invariant for the group chain $\cG$ depends on the section, so can be considered as indeterminately defined. Thus, the usage in Definition~\ref{def-asympdisc2}  is   more precise.

 The discussion above shows the following:
 \begin{prop}\label{prop-ADequivalence}
 Let  $(X_{\infty} , G_0 , \Phi_0)$   be the equicontinuous    minimal Cantor system associated to the group chain $\cG  = \{G_{\ell}\}_{\ell \geq 0}$. Then the asymptotic class $[\cI(\Phi_0)]_{\infty}$ for the sequence of surjective homomorphisms defined by \eqref{eq-surjectiveUn} equals the asymptotic class $[\cD(\cG)]_{\infty}$ for the sequence of surjective homomorphisms  $\cD(\cG)$ defined by \eqref{eq-surjectiveDn}.
 \end{prop}
 
  It remains to show that there exists a plethora of examples of wild group chains for which the asymptotic discriminant invariants are distinct, so that this is an effective invariant of weak solenoids. We construct such examples   in the next two sections.

\section{General constructions of solenoids}\label{sec-construction}

In this section, we recall two general constructions which are used to produce examples of weak solenoids with prescribed properties. The first is a well-known construction which yields a foliated space  from a group action, called the \emph{suspension construction}, as discussed in   \cite[Chapter 3]{CandelConlon2000} for example, or Chapter 6 of \cite{CN1985} in the case of smooth foliations. 
In Section~\ref{subsec-suspension} we describe this construction in the   context of   Cantor actions. 

The second construction produces group chains whose associated equicontinuous minimal Cantor systems  $(X_{\infty} , G_0 , \Phi_0)$  have prescribed discriminant invariants. This is described in Section~\ref{subsec-lenstra}, and is inspired by a construction   introduced by Fokkink and Oversteegen \cite[Lemma~37]{FO2002} and attributed to Hendrik Lenstra. This method was  used in \cite[Section~10]{DHL2016c} to construct examples of stable but not homogeneous solenoids. In   Section~\ref{sec-wildex},  the methods of Section~\ref{subsec-lenstra}  are used  to construct examples of wild actions, and thus by the methods of Section~\ref{subsec-suspension} we obtain examples of wild solenoids.

\subsection{The suspension construction}\label{subsec-suspension}

Let $X$ be a Cantor space, and $H$ a \emph{finitely-generated} group, though not necessarily finitely-presented, and assume there is given  an action $\vp \colon H \to Homeo(X)$. Suppose that $H$ admits a generating set $\{g_1 , \ldots , g_k\}$, let $\alpha_k \colon \mZ * \cdots * \mZ \twoheadrightarrow H$ be the surjective homomorphism from the free group on $k$ generators to $H$, given by mapping generators to generators.  
 Next, let $\Sigma_k$ be a compact surface without boundary of genus $k$. For a choice of basepoint $x_0 \in \Sigma_k$, set   $G_0 = \pi_1(\Sigma_k, x_0)$,  then there is a   homomorphism $\beta_k  \colon G_0 \to \mZ * \cdots * \mZ$ onto the free group of $k$ generators.   
The composition of these   maps  yields a   homomorphism 
\begin{equation}\label{eq-suspensionmap}
\Phi_0 = \vp \circ \alpha_k \circ \beta_k  \colon G_0 = \pi_1(\Sigma_k, x_0) \to  \mZ * \cdots * \mZ \to H \to  Homeo(X) \ .
\end{equation}
 
Next, let $\wtSigma_k$ denote the universal covering space of $\Sigma_k$, equipped with the right action of $G_0$ by covering transformations. Form the product space $\wtSigma_k \times X$ which has a foliation $\wtF$ whose leaves are   the slices $\wtSigma_k \times \{x\}$ for each $x \in X$. Define a left action of $G_0$ on $\wtSigma_k \times X$, which for $g \in G_0$ is given  by $(y,x) \cdot g = (y \cdot g, \Phi_0(g^{-1})(x))$. For each $g$,  this action preserves the foliation $\wtF$, so we obtain a foliation $\FfM$ on the quotient space $\fM = (\wtSigma_k \times X)/G_0$.  Note that all leaves of $\FfM$ are surfaces, which are in general non-compact.

The space $\fM$ is a foliated Cantor bundle over $\Sigma_k$, and the holonomy of this bundle $\pi \colon \fM \to \Sigma_k$ acting on the fiber $V_0 = \pi^{-1}(x_0)$ is canonically identified with the action $\Phi_0 \colon G_0 \to Homeo(X)$. Consequently,  if the action $\Phi_0$ is minimal  then the foliation $\FfM$ is minimal. If the action $\Phi_0$ is equicontinuous   then $\FfM$ is an equicontinuous foliation in the sense of Definition~\ref{def-equicontinuous}.
 
 There is a variation of the above construction, where we assume that $G_0$ is a   \emph{finitely-presented} group, and there is given a homomorphism   $\Phi_0 \colon G_0 \to Homeo(X)$.
 Then there is a   folklore result (for example, see \cite{Massey1991}) that there   exists a closed connected $4$-manifold $B$ such that for a choice of basepoint $b_0 \in B$, the fundamental group $\pi_1(B, b_0)$ is isomorphic to $G_0$.  Then the suspension construction can be applied to the homomorphism $\Phi_0 \colon \pi_1(B, b_0) \to Homeo(X)$, where we replace $\Sigma_k$ above with $B$,  and the space $\wtSigma_k$ with the universal covering $\widetilde{B}$ of $B$.

\subsection{The Lenstra construction}\label{subsec-lenstra} We give a general   construction of a group chain with prescribed core and discriminant groups.  See \cite[Section~10]{DHL2016c} for further discussion of this construction, and the proof of Lemma~37 in \cite{FO2002}  for the motivation behind this approach to constructing group chains.

Let  $\whH$ be a given profinite group, with  a finitely-generated dense subgroup $H  \subset  \whH$.
Let $\cD \subset \whH$ be a closed subgroup of infinite index, and so $\cD$ is either a finite group or is a profinite group.  

The quotient       $X = \whH/\cD$ is a Cantor space.
As $H$ is dense in $\whH$, the left action of $H$ on $X$ is minimal. 
The \emph{rational  core} in $\whH$ of $\cD$ is defined to be the subgroup   
\begin{equation}\label{eq-rationalcorepro}
{\rm core}_{H} \  \cD =    \bigcap_{h \in H} h   \cD   h^{-1} \subset \whH \ .
\end{equation}
 
\begin{lemma}\label{lem-rationalcore}
The left action of $\whH$ on   $X$ is effective if and only if ${\rm core}_{H} \ \cD = \{\whe\}$.
\end{lemma}
\proof
Suppose there exists  $\whe \ne \whk \in {\rm core}_{H} \ \cD$, then $h^{-1}  \whk h     \in   \cD$ for all $h \in H$, and so $\whk \cdot h \cD  = h \cD$ for all $h \in H$. Thus, the continuous action of $\whk$ fixes the dense subset  $\{h \cD \mid h \in H\} \subset X$, hence must act trivially on all of $X$. That is, the action of $\whH$ on $X$ is not effective. Thus, if the action is effective, then ${\rm core}_{H} \ \cD   = \{\whe\}$. 

To show the converse, assume that ${\rm core}_{H} \ \cD   = \{\whe\}$,  then observe that 
$$ {\rm core}_{\whH} \ \cD =  \bigcap_{\whh \in \whH} \ \whh   \cD   \whh^{-1} \subset  \bigcap_{h \in H} h    \cD   h^{-1}   =  {\rm core}_{H} \ \cD = \{\whe\} \ .$$
Thus,  the intersection of all isotropy groups $\{\whh   \cD   \whh^{-1}  \mid \whh \in \whH\}$  for the action of $\whH$ on $X$ is trivial, hence the action is effective.
 \endproof

We  next construct   a group chain $\cH = \{H_{\ell}\}_{\ell \geq 0}$  with $H_0 = H$ and discriminant group  isomorphic to    $\cD$.
By the assumption that $\whH$ is a profinite group, there exists a  clopen neighborhood system  $\{\whU_{\ell}  \mid \ell \geq 1\}$  about the identity in $\whH$, so that: 
\begin{enumerate}
\item each $\whU_{\ell}$ is a normal subgroup of $\whH$ ;
\item for each $\ell \geq 0$ there is a proper inclusion $\whU_{\ell+1}  \subset \whU_{\ell}$ ;
\item the  intersection ${\ds \bigcap_{\ell \geq 1} \ \whU_{\ell}   = \{\whe\}}$ .
\end{enumerate}
 In particular, each quotient $\whH/\whU_{\ell}$ is a finite group. Let $\iota^{\ell+1}_{\ell} \colon \whH/\whU_{\ell +1} \to \whH/\whU_{\ell}$ be the map induced by the inclusions of cosets.
 The sequence of quotient maps 
 \begin{align}\label{eq-quotientul}\{q_{\ell} \colon \whH \to \whH/\whU_{\ell} \mid \ell \geq 1\}\end{align}  
 induces an isomorphism $\iota_{\infty}$ of $\whH$ with the inverse limit group 
 \begin{equation}
\iota_{\infty} \colon \whH \cong H_{\infty} \equiv \lim_{\longleftarrow}\, \left \{\iota^{\ell+1}_{\ell} \colon \whH/\whU_{\ell +1} \to \whH/\whU_{\ell}  \right\} \ .
\end{equation}
Next, for each $\ell \geq 1$, set $\whW_{\ell} =  \whU_{\ell} \cdot \cD$ which is a   subgroup of $\whH$, as $\whU_{\ell}$ is a normal subgroup. 
 Moreover, the assumption that $\cD$ is compact implies that each $\whW_{\ell}$ is a clopen subset of $\whH$. Set  $H_{\ell} = H \cap \whW_{\ell}$ which is a subgroup of finite index in $H$, and so $\cH = \{H_{\ell}\}_{\ell  \geq 0}$ is a group chain in $H$.

Next, for each $\ell \geq 0$,  set $X_{\ell} = H/H_{\ell}$ which is a finite set with a transitive left action of $H$. There are quotient maps $p_{\ell +1} \colon X_{\ell +1} \to X_{\ell}$ which commute with the $G$-action, and introduce the inverse limit space
  \begin{equation}\label{eq-algmodel}
X_{\infty} = \lim_{\longleftarrow} ~ \{ p_{\ell +1} \colon X_{\ell +1}  \to X_{\ell}\} ~ \subset \prod_{\ell \geq 0} ~ X_{\ell}  \ ,
\end{equation}
which  is   a totally disconnected compact perfect set, so is a Cantor set. The diagonal left action of $H$ on $X_{\infty}$ is   equicontinuous and minimal, so we obtain an equicontinuous minimal Cantor system $(X_{\infty} , H, \Phi)$. 
Note that  $H$ is dense in $\whH$ implies that 
\begin{align}\label{eq-equalquotients}X_{\ell} = H/H_{\ell} \cong \whH/\whW_{\ell}.\end{align}
Thus,  
 there are surjections
 \begin{equation}\label{eq-surjectionstau}
\tau_{\ell} \colon X \equiv \whH/\cD \to    \whH  /\whW_{\ell} \cong    H / H_{\ell}  = X_{\ell}
\end{equation}
 such that $p_{\ell +1} \circ \tau_{\ell +1} = \tau_{\ell}$. The collection of maps $\{\tau_{\ell} \colon X \to   X_{\ell} \mid \ell \geq 1\}$
   then induces an $H$-equivariant homeomorphism $\tau_{\infty} \colon X  \cong X_{\infty}$.
In particular,  for each $y = (g_{\ell} H_{\ell}) \in X_{\infty}$ there is a   sequence $(h_{\ell}) \in \whH$ such that for the sequence of cosets $(h_{\ell} \cD) \in \whH/\cD$ we have $\tau_{\infty}(h_{\ell} \cD) = (g_{\ell} H_{\ell})$.  

 We next relate the closed group $\cD \subset \whH$ to the discriminant group $\cD_{\infty}$ for the group chain   $\cH = \{H_{\ell} \}_{\ell \geq 0}$.   
Recall that $\ds  C_{\ell}   = {\rm core}_{H} \  H_{\ell}   =    \bigcap_{g \in H} \ g H_{\ell} g^{-1} \subset H$  
is the   core in $H$ of  $H_{\ell}$, 
and  
\begin{equation}\label{eq-discpro}
\cD_{\infty} \equiv \lim_{\longleftarrow}\, \left \{\delta^{\ell+1}_{\ell} \colon  H_{\ell+1}/C_{\ell+1} \to H_{\ell}/C_{\ell}  \right\}   \   .
\end{equation}

\begin{lemma} 
There is an isomorphism $\cD_\infty \to \cD$. 
\end{lemma}

\proof Consider the quotient maps \eqref{eq-quotientul}, that is, for $\ell \geq 1$, we have $q_{\ell} \colon \whH \to  \whH/\whU_{\ell}$, where  $\whH/\whU_{\ell}$ is a finite group.   Since $H$ is dense in $\whH$, it has non-trivial intersection with  each clopen set $g  \whU_{\ell}$. Then 
\begin{equation}\label{eq-discell}
\cD^{\ell} \equiv q_{\ell}(\cD) = q_{\ell}(\whW_{\ell}) = q_{\ell}(H \cap \whW_{\ell}) = q_{\ell}(H_{\ell})  \ ,
\end{equation}
where the first equality holds because $\whW_\ell = \whU_\ell \cdot \cD$. As $q_\ell$ are homomorphisms, and the rational core of $\cD$ is trivial by Lemma~\ref{lem-rationalcore},    we   have
\begin{eqnarray}\label{eq-coreell}
\lefteqn{q_{\ell}(C_{\ell}) = q_{\ell}(\bigcap_{g \in H} \ g \ H_{\ell} \ g^{-1}) 
  =   \bigcap_{g \in H} \ q_{\ell}(g) \ q_{\ell}(H_{\ell}) \ q_{\ell}(g)^{-1}  }\\ & = &  \bigcap_{g \in H} \ q_{\ell}(g) \ q_{\ell}(\cD) \ q_{\ell}(g)^{-1} = q_\ell( \bigcap_{g \in H} g\cD g^{-1})
   = \{ e_{\ell}  \}  \in \whH/\whU_{\ell} \ . \nonumber
\end{eqnarray}
Since by definition $C_\ell$ is a subgroup of $H_\ell \subset H$, it follows that $C_{\ell} \subseteq H \cap \whU_{\ell}$. The subgroup $\whU_{\ell}$ is normal in $\whH$, so $H \cap \whU_{\ell}$ is normal in $H_{\ell} = H \cap \whW_{\ell}$, hence $H \cap \whU_\ell \subseteq C_{\ell}$, as $C_\ell$ is maximal. We then obtain $C_{\ell} = H \cap \whU_{\ell} = H_\ell \cap \whU_\ell$. 
 In particular,  we obtain induced maps on the quotients, $\ovq_{\ell} \colon H/C_{\ell} \to \whH/\whU_{\ell}$, and so by \eqref{eq-discell} we have that 
 $\ovq_{\ell}(H_{\ell}/C_{\ell}) = q_{\ell}(\cD) = \cD^{\ell}$ for all $\ell \geq 0$.  

The map $\iota^{\ell+1}_{\ell} \colon H_{\ell+1} \to H_{\ell}$  induces a map  $\iota^{\ell+1}_{\ell} \colon \cD^{\ell+1} \to \cD^{\ell}$, and so 
  for the inverse limit  we have
\begin{equation}\label{eq-discpro2}
\cD_{\infty} =  \lim_{\longleftarrow}\, \left \{\delta^{\ell+1}_{\ell} \colon  H_{\ell+1}/C_{\ell+1} \to H_{\ell}/C_{\ell}  \right\}  \cong \lim_{\longleftarrow}\, \left \{\iota^{\ell+1}_{\ell} \colon  \cD^{\ell+1} \to \cD^{\ell}   \right\}  \cong \cD \   , 
\end{equation}
where the isomorphism on the right hand side of \eqref{eq-discpro2}  follows as $\{\whU_{\ell} \mid \ell \geq 1\}$   is a clopen neighborhood system about the identity in $\whH$. 
 \endproof

 \subsection{Wild criteria}\label{subsec-wildcriteria} 
 We next formulate the stability property in Definition~\ref{def-stableGC} in terms  of the subgroup $\cD$, and obtain   criteria which are sufficient to imply that  the   
 equicontinuous minimal Cantor system $(X_{\infty} , H, \Phi)$ is wild. 
 First, note that for   $n \geq 1$, we have $\cD \subset \whW_{n} = \cD \cdot \whU_{n} =   \whU_{n} \cdot \cD$, and define 
 $U_{n} \equiv  \whW_{n}/\cD   \subset \whH/\cD = X$, 
which  is a clopen neighborhood of the basepoint $x_0 = \whe \cD \in   X$. 
The proof of the following is immediate:
\begin{lemma}\label{lem-isotroptyn}
For   $\ell \geq 1$,  $g \in H$ satisfies $g \cdot U_{\ell} = U_{\ell}$ if and only if $g \in H_{\ell}$.
\end{lemma}

Lemma~\ref{lem-isotroptyn} implies that we can identify $U_n$ with the inverse limit in \eqref{Ginfty-definen}. We now consider the Ellis group and the discriminant group of the truncated group chain 
$\cH_n = \{H_{\ell}\}_{\ell \geq n}$ as in \eqref{eq-discquotients1}.
  For each $n \geq 1$, let $\rho_n \colon \whW_n \to Aut(\whW_{n})$ which for $\whh \in \whW_n$ and $\whg \in \whW_{n}$ is defined by 
$\rho_n(\whh)(\whg) = \whh \ \whg \ \whh^{-1}$, and let 
$\whrho_n = \rho_n|\cD \colon \cD \to Aut(\whW_{n})$ be it's restriction to $\cD$.

\begin{lemma}\label{lem-trivialaction}
Let $\whh \in \cD$ then $\rho_n(\whh) \in Aut(\whW_{n})$ is the identity implies that the left action of $\whh$  on $U_n$ is the identity. 
\end{lemma}
\proof Assume that  $\rho_n(\whh)$ is the identity, then  $\whh  \ \whg \ \whh^{-1} = \whg$ for every $\whg \in \whW_n$. Then for the action on cosets of $\cD$ we have $\whh \cdot \whg \cD = \whh \whg \cD = \whg \whh \cD = \whg \cD$ as $\whh \in \cD$. 
\endproof

Let $\cK_n = ker(\whrho_n) \subset \cD$ denote the kernel of the representation $\whrho_n$. Then $m > n$ implies   $\cK_n \subseteq \cK_m$.

We would like to obtain   representations  of the Ellis group of the restricted action on $U_n$  and it's discriminant group, in terms of $\whW_n$, $\cD$ and $\cK_n$.
Note that the converse of Lemma \ref{lem-trivialaction} does not hold, as the identity 
  $$ \whd \cdot \whg \cD = \whd \whg  \whd^{-1}\cD = \whg \cD $$
only implies that the adjoint action of $\whd$   permutes points within cosets of $\cD$ in $\whW_n$ but need not act as the identity   on $\cD$. Thus the group $\cD/\cK_n$ may act on $U_n$ non-effectively and cannot be identified with the  isotropy group of the Ellis group action, and in particular, this means that $\cD/\cK_n$ may have non-trivial core in $\whW_n/\cK_n$.

We want to find a condition which eliminates this situation. Consider the sequence 
$$\whW_0 = \whH \supset \whW_1 \supset \cdots \supset \whW_n \supset \cdots \ .$$ 
The discriminant group $\cD$ of the action of $\whW_0$ on $X = \whW_0/\cD$ has trivial core by Proposition~\ref{prop-discrisnotnormal}. 
Consider the action of $\whW_n$ on $U_n = \whW_n/\cD$. Note that   the action restricted to $\cD$ induces an action of   $\cD/\cK_n$ by Lemma~\ref{lem-trivialaction}.  
Suppose there exists $\whg \in \whW_n$ such that the action of $\whg$ on $U_n$ is trivial, but $\whg \notin \cK_n$. 
Then for each $\whh \in \whW_n$ we have $\whg \whh \cD = \whh \cD$, so $\whg \whh \cD \whh^{-1} = \whh \cD \whh^{-1} $, that is, $\whg \in \bigcap_{\whh \in \whW_n} \whh \cD \whh^{-1} $. 
Thus,  the discriminant group of the restricted action is isomorphic to $\cD/\cK_n$ if the following condition is satisfied:
  \begin{equation}\label{non-triv-condition}
  {\rm core}_{\whW_n} \ \cD/\cK_n =\bigcap_{\whh \in \whW_n} \whh \, (\cD/\cK_n)\ \whh^{-1}  = \{\whe_n\} \ .
  \end{equation}

Recall that there are quotient maps $q_n \colon \whH \to \whH/\whU_{n}$, and denote $H_{n,n} = q_n(\whW_n) = q_n(H_n) = \cD^n$, 
where the equalities follow from  \eqref{eq-discell}. 
Consider the truncated chain $\cH_n = \{H_\ell\}_{\ell \geq n}$, and denote $H_{n,\ell} = q_\ell(H_n)$. 
Also, we have $\cD^\ell = q_\ell (\cD) = q_\ell(H_\ell) $. A computation similar to \eqref{eq-coreell} shows that
\begin{equation*}
  q_{\ell}(\bigcap_{g \in H_n} \ g \ H_{\ell} \ g^{-1}) 
  =   \bigcap_{g \in H_{n,\ell}} \ k \ \cD^\ell\ k^{-1}  = {\rm core}_{H_{n,\ell}} \ \cD^\ell \ ,  
  \end{equation*}
and then \eqref{non-triv-condition} implies that for all $\ell > n >1$ we must have
   \begin{equation}\label{eq-corenm}
{\rm core}_{H_{n, \ell}} \  \cD^\ell =    \bigcap_{k \in H_{n, \ell}} k   \cD^\ell   k^{-1} = \{e_\ell\} \in H_{n, \ell}   \ .
\end{equation}

We  then have the following sufficient criteria for  an action to be wild:

\begin{prop}\label{prop-wildcriterium}
Suppose that the $H_{n, \ell}$ core of $\cD^\ell$ is trivial in $\whW_{\ell}/\whU_\ell$ for all $\ell > n \geq 1$, and that 
 the chain  $\{\cK_n \mid n \geq 1\}$ is not bounded above. Then the equicontinuous minimal Cantor system $(X_{\infty} , G, \Phi)$ is wild.
\end{prop}
\proof
By assumption, there exists an increasing sequence of integers $\{n_{\ell} \mid \ell \geq 1\}$ such that  for all $\ell \geq 1$,  the inclusion $\cK_{n_{\ell}} \subset \cK_{n_{\ell +1}}$ is  proper. Thus, there exists a sequence of elements $\whd_{\ell+1} \in \whW_{n_{\ell}}$ such that $\whd_{\ell+1} \in \cK_{n_{\ell +1}}$ but $\whd_{\ell+1} \not\in \cK_{n_{\ell}}$. The assumption on the triviality of $H_{n, \ell}$-cores implies that the discriminant group of the restricted action is isomorphic to $\cD/\cK_n$. By assumption on the kernels $\cK_n$ the maps $\psi_{n,n+1}: \cD/\cK_n \to \cD/\cK_{n+1}$ of the discriminant groups, defined in Section \ref{subsec-ADGC}, are not injective for any $n \geq 0$, and it follows that the action of $G$ on $X$ is wild. 
\endproof

\begin{remark}\label{embeddings}
{\rm We note that the actions satisfying the conditions of Proposition \ref{prop-wildcriterium} are not LCSQA. Indeed, consider the collection of elements $\{\whd_{\ell}\}$ chosen as in the proof of the Proposition. Since $\whd_{\ell+1} \in \cK_{n_{\ell+1}}$, by Lemma \ref{lem-trivialaction} $\whd_{\ell+1}$ acts trivially on the clopen set $U_{n_{\ell+1}} = \whW_{\ell+1}/\cD$, and non-trivially on the larger set $U_{n_{\ell}} = \whW_{\ell}/\cD$. Since $\bigcap U_{n_\ell}$ is a singleton, the action is not LCSQA.
}
\end{remark}

We conclude with some remarks about wild actions. Recall that 
the kernel for the group chain $\cH$ is     
\begin{equation}\label{eq-kernelpro}
K(\cH) = \bigcap_{\ell \geq 0} \  H_{\ell} =  \bigcap_{\ell \geq 0} \ \left(  H \cap \whW_{\ell} \right) = H \ \cap  \   \bigcap_{\ell \geq 0} \    \whW_{\ell}   = H \cap \cD   \ . 
\end{equation}

\begin{lemma}
Let $\cD \subset \whH$ be a closed subgroup, and 
let  $\cH = \{H_{\ell}\}_{\ell \geq 0}$ with $H_0 = H \subset \whH$ be the group chain associated with the choice of a clopen neighborhood system $\{U_{\ell}  \mid \ell \geq 1\}$ about the identity in $\whH$. 
Given $y = (g_{\ell} H_{\ell}) \in X = \whH/\cD$  let $\cH^y = \{g_{\ell} H_{\ell} g_{\ell}^{-1}\}_{\ell \geq 0}$ be the conjugate group chain to $\cH$. Then for $\whh = (h_{\ell}) \in \whH$ such that   $\tau(\whh \cD) = y$, we have
\begin{equation}
K(\cH^y) =  H  \ \cap \ (\whh \ \cD \ \whh^{-1})  \ .
\end{equation}
 In particular, if $\cD$   satisfies $H \cap \whg \cD \whg^{-1} = \{e\}$ for all $\whg \in \whH$, then the kernel of each group chain $\cH^y$ centered at $y \in X$ is trivial. 
\end{lemma}
 
 If the kernel $K(\cH^y)$ is trivial for all $y \in X$ then the action of $H$ on $X$ satisfies the SQA property.  It may still happen that the action of the Ellis group $\whH$ of the action is not LCSQA, the hypotheses of Proposition~\ref{prop-wildcriterium} are satisfied, and so the action of $H$ on $X$ is wild.

\section{Non-homeomorphic  wild solenoids}\label{sec-wildex}

  In this section, we give a new construction  of weak solenoids whose monodromy actions are wild,  inspired by the constructions of stable solenoids in  Section~10.5 of \cite{DHL2016c}.  More precisely, we construct  
  equicontinuous minimal Cantor actions which are wild, then use the suspension construction described in  Section~\ref{subsec-suspension} to realize these actions as the global monodromy actions of weak solenoids. 
 
 The following examples are inspired by the constructions of Lubotzky in \cite{Lubotzky1993} of compact   torsion subgroups in the profinite completion of  finitely generated, torsion free groups. We do not require the full strength of Lubotzky's methods,  though the reader may find the discussion in Section~10.5 of \cite{DHL2016c} helpful for establishing the context. 
 
  \subsection{Wild actions}\label{subsec-Lwild}
  We first     recall some basic facts as used in \cite{Lubotzky1993}.   For $N \geq 3$, let $\G_N = {\bf SL}_N(\mZ)$ denote the $N \times N$ matrices with integer entries.
       The group $\G_N$ is finitely-generated and residually finite, and hence so are all finite index subgroups of $\G_N$.  For an integer $M \geq 2$, 
let $\G_N(M)$ denote the congruence subgroup 
    $$\G_N(M) \equiv  {\rm Ker} \left\{\varphi_M \colon  {\bf SL}_M(\mZ) \to {\bf SL}_N(\mZ/M \mZ) \right\} \ .$$
For $M \geq 3$, $\G_N(M)$ is torsion-free.  
 Moreover,  by the congruence subgroup property, every finite index subgroup of $\G_N$ contains $\G_N(M)$ for some non-zero $M$;  that is, congruence subgroups form a neighborhood basis about the identity for the profinite topology of ${\bf SL}_N(\mZ)$.    
 This implies that
 \begin{equation}\label{eq-prodprimes}
\widehat{{\bf SL}_N(\mZ)} \equiv \lim_{\longleftarrow} ~ {\bf SL}_N(\mZ/M \mZ) \cong {\bf SL}_N(\widehat{\mZ}) \cong  \prod_{p \in \cP} \ {\bf SL}_N(\widehat{\mZ}_p)
\end{equation}
where $\widehat{\mZ}$ is the profinite completion of $\mZ$, and $\widehat{\mZ}_p$ is the   $p$-adic completion of $\mZ$. Here,      $\cP = \{p_i \mid i \geq 1\}$ is the set of primes listed in increasing order, where $p_1 = 2$ and $p_2 =3$.  The  second isomorphism  uses the observation that   $\ds \widehat{\mZ} = \prod_{p \in \cP}  \ \widehat{\mZ}_p$. 
Note that the matrix-group factors in the Cartesian product on the right hand side of  \eqref{eq-prodprimes} commute with each other.
   
Fix an integer $N \geq 3$.   Let $G \subset \G_N$ be a   finite index, torsion-free subgroup, which is then finitely-generated, 
  and its profinite completion $\whG$ is a clopen subgroup of finite index in $\widehat{{\bf SL}_N(\mZ)}$. 
  Then there exists an index $k_G \geq 2$ such that for $\cP(G) = \{p_i \mid i \geq k_G\}$  we have 
 \begin{equation}\label{eq-embedp}
 \prod_{i \geq k_G} \ {\bf SL}_N(\widehat{\mZ}_{p_i}) \subset \whG \ ,
\end{equation}
where the inclusion is via the isomorphism in \eqref{eq-prodprimes}.  In particular, this means that the diagonal embedding of $G$ is dense in the product of groups on the left-hand-side in \eqref{eq-embedp}.
   
  Let $\ds \whH = \prod_{i \geq k_G} \ {\bf SL}_N(\mZ/p_i\mZ)$ be the infinite product of finite groups with the Tychonoff product topology, which is thus a Cantor group. The group ${\bf SL}_N(\mZ)$ maps onto   each factor ${\bf SL}_N(\mZ/p_i\mZ)$, and let $Q \colon {\bf SL}_N(\mZ) \to  \whH$ denote the product of all the maps. 
  
  Let $H = Q(G) \subset \whH$ denote the image of $G$. Then $H$ is a finitely-generated subgroup of $\whH$.

  \begin{lemma}
  The subgroup $H \subset \whH$ is dense in $\whH$.
  \end{lemma}
  \proof For $p_i \in \cP(G)$ consider the projections 
  \begin{equation}
\kappa_i \colon \widehat{{\bf SL}_N(\mZ)} \stackrel{\theta_i}{\longrightarrow} {\bf SL}_N(\widehat{\mZ}_{p_i})  \stackrel{\alpha_i}{\longrightarrow} {\bf SL}_N(\mZ/p_i \mZ)
\end{equation}
where   $\theta_i$ is defined via the isomorphism \eqref{eq-prodprimes},  $\alpha_i$ is induced by the quotient map  $\widehat{\mZ}_{p_i} \to \mZ/p_i \mZ$, and  $\kappa_i = \alpha_i \circ \theta_i$. 
    Since $G \subset {\bf SL}_N({\mZ}) $ is dense in the product  \eqref{eq-embedp}, the restriction $\theta_i|G$  has dense image, 
    and then $\kappa_i| G = \alpha_i \circ \theta_i|G$ is surjective.    
    Also, $\kappa_i$ factors through the projection ${\rm pr}_i :  \prod_{p_i \in \cP(G)} \ {\bf SL}_N(\widehat{\mZ}_{p_i}) \to {\bf SL}_N(\mZ/p_i \mZ) $, and so the restriction ${\rm pr}_i |G$ is surjective.
  
The embeddings $Q: G \to \whH$ and $G \to \prod_{p_i \in \cP} \ {\bf SL}_N(\widehat{\mZ}_{p_i})$ are diagonal, therefore, the following diagram commutes
   \begin{equation*} 
   \xymatrix{ G \ar[r] \ar[d]_Q & \prod_{p_i \in \cP(G)} \ {\bf SL}_N(\widehat{\mZ}_{p_i}) \ar[r]^{{\rm pr}_i} &  {\bf SL}_N(\mZ/p_i\mZ) \ar[d]^{\rm id}  \\ 
   \whH     \ar[r]            & \prod_{p_i \in \cP(G)} \ {\bf SL}_N(\mZ/p_i\mZ) \ar[r]                                  & {\bf SL}_N(\mZ/p_i\mZ) } 
   \end{equation*}
   and so the projection of $H = Q(G)$ onto $ {\bf SL}_N(\mZ/p_i\mZ)$ is surjective. It follows that $H$ is dense in $\whH$ for the product Tychonoff topology.
\endproof
  
    Note that  the kernel of the restricted map $Q \colon G \to \whH$ is an infinitely-generated subgroup of $G$.

  For each $p_i \in \cP(G)$, choose a non-trivial subgroup $A_{p_i}  \subset {\bf SL}_N(\mZ/p_i\mZ)$ with trivial core. 
  That is, $A_{p_i}$ is chosen so that the   trivial group is the largest normal subgroup of ${\bf SL}_N(\mZ/p_i\mZ)$ contained in $A_{p_i}$.  For example, the alternating group on $N$ symbols embeds into $ {\bf SL}_N(\mZ)$, and so it embeds into $ {\bf SL}_N(\mZ/p_i \mZ)$. For $N=3$ or $N \geq 5$, the alternating group on $N$ symbols is simple, so $A_{p_i}$ maybe chosen to be any non-trivial  proper subgroup of the   image of the alternating group in ${\bf SL}_N(\mZ/p_i \mZ)$.
  
   Since $A_{p_i}$ has trivial core, it follows that each $g \in {\bf SL}_N(\mZ/p_i\mZ)$ with $g \notin A_{p_i}$ acts non-trivially on the left on $Y_{p_i} =  {\bf SL}_N(\mZ/p_i\mZ)/A_{p_i}$.  Then define
  $$\cD = \prod_{i \geq k_G} \ A_{p_i} \subset \whH \ ,$$
and introduce the quotient space  
  $$X = \whH/\cD  = \prod_{i \geq k_G} \ {\bf SL}_N(\mZ/p_i\mZ)/A_{p_i}  \cong   \prod_{i \geq k_G}  \ Y_{p_i} \ .$$
Define the   action of $G$ on $X$ by setting, for $g \in G$ and $\whx \in \whH/\cD$, 
$\Phi_0(g)(\whx) = Q(g) \cdot \whx$. Then we obtain
an  equicontinuous minimal Cantor action $(X, G ,  \Phi_0)$. Note that the kernel of the homomorphism $\Phi_0 \colon G \to Homeo(X)$ is an infinitely generated subgroup.

 We next calculate the discriminant groups for $(X, G ,  \Phi_0)$, using   the criteria given in Proposition~\ref{prop-wildcriterium}.  
 
For $\ell > k_G$  define 
\begin{equation}\label{eq-Uell}
\whU_{\ell} = \prod_{k_G \leq i < \ell} \ \{e_{p_i}\} ~ \times ~ \prod_{i \geq \ell} \ {\bf SL}_N(\mZ/p_i\mZ) \quad ; \quad U_{\ell} = \prod_{k_G \leq i < \ell} \ \{e_{p_i}\} ~ \times ~ \prod_{i \geq \ell} \ Y_{p_i} \ ,
\end{equation}
 where  $e_{p_i} \in A_{p_i} \subset {\bf SL}_N(\mZ/{p_i}\mZ)$ is the identity element. 
 Each $\whU_{\ell}$ is a clopen subset of $\whH$, and $U_{\ell}$ is a clopen neighborhood of $x_0 \in X$.
 Thus, we obtain a   descending chain of clopen neighborhoods $\{\whU_{\ell} \mid \ell > k_G \}$ of the identity in $\whH$. 
Observe that for $\ell > k_G$ 
 $$\whW_{\ell} =  \whU_{\ell} \cdot \cD = \prod_{k_G \leq i < \ell} \ A_{p_i} ~ \times ~ \prod_{i \geq \ell} \ {\bf SL}_N(\mZ/p_i\mZ)\ .$$
Set    
$\ds H_{\ell}  = H \cap \whW_{\ell}$, so that   by \eqref{eq-equalquotients}  the finite set $X_{\ell}$ has the form
$$X_{\ell} = \whH/\whW_{\ell} \cong \prod_{k_G \leq i < \ell} \ {\bf SL}_N(\mZ/p_i \mZ)/A_{p_i}  \cong \prod_{k_G \leq i < \ell} \ Y_{p_i}  \ . $$
Observe that   each $\ell \geq k_G$, the kernel $\cK_{\ell}$ of the adjoint representation  $\whrho_{\ell} \colon \cD \to Aut(\whW_{\ell})$ is given by 
\begin{align}\label{kernel-compute}\cK_{\ell} = \prod_{k_G \leq i < \ell} \ A_{p_i} ~ \times ~ \prod_{i \geq \ell} \ \{e_{p_i}\} \ ,\end{align}
so that $\{\cK_{\ell} \mid \ell > k_G \} \subset \cD$ form an increasing chain of subgroups of $\cD$ whose union is dense in $\cD$.

Next, for $\ell > k_G$, we have the projection onto the quotient group,   as in \eqref{eq-quotientul},
 $$q_{\ell} \colon  \whH \to \whH/\whU_{\ell}   \cong \prod_{k_G \leq i < \ell} \ {\bf SL}_N(\mZ/p_i \mZ)  \ .$$

For $\ell > n > k_G$, recall that $H_{n}  = H \cap \whW_{n}$ and      $\whW_{n} =  \whU_{n} \cdot \cD$ where $\whU_n$ is defined by \eqref{eq-Uell}, then set  
$$H_{n, \ell} = q_\ell(H_{n}) = q_\ell(\whW_n) \cong  \prod_{k_G \leq i < n} \ A_{p_i} ~ \times ~  \prod_{n \leq i < \ell} \ {\bf SL}_N(\mZ/p_i \mZ) \subset \prod_{k_G \leq i < \ell} \ {\bf SL}_N(\mZ/p_i \mZ)  \ .$$
Finally, the image under the quotient map $q_n \colon \whH \to \whH/\whU_{n}$  of the group $\cD$ is given by 
$$ q_\ell(\cD) =  \prod_{k_G \leq i < \ell} \ A_{p_i}    \subset \prod_{k_G \leq i < \ell} \ {\bf SL}_N(\mZ/p_i \mZ) \ ,$$
and by \eqref{kernel-compute}
 \begin{align}\label{discr-restricted} q_\ell(\cD/\cK_n) =  \prod_{n \leq i < \ell} \ A_{p_i}    \subset \prod_{k_G \leq i < \ell} \ {\bf SL}_N(\mZ/p_i \mZ) \ . \end{align}

By the choice of each subgroup $A_{p_i}$ we have that the    $H_{n, \ell}$ core of $q_\ell(\cD/\cK_n)$  is trivial, and so by Proposition~\ref{prop-wildcriterium} the discriminant group $\cD_{n,n}$ of the action, restricted to $U_n = \whW_n/\cD$, is isomorphic to $\cD/\cK_n$. Then, taking the inverse limit of groups in \eqref{discr-restricted}, we obtain that
\begin{equation}\label{restricted-discr}
\cD_{n} = \prod_{i \geq n} \ A_{p_i} \subset \whH \ ,
\end{equation}
and the surjective homomorphisms $\psi_{n} \colon \cD_{n} \to \cD_{n+1}$ in  \eqref{eq-surjectiveDn} are given by 
\begin{equation}\label{eq-surjectiveDDn}
\psi_{n,n+1} \colon \prod_{i \geq n} \ A_{p_i} \to \prod_{i \geq n+1} \ A_{p_i}
\end{equation}
In particular, we have $ker(\psi_{n,n+1}) = A_{p_n}$ for all $n \geq 1$. We have thus shown:

\begin{prop}\label{prop-ex1wild}
The action  $(X, G, \Phi_0)$     is wild, with   asymptotic discriminant given by the asymptotic class of the sequence of surjective homomorphisms 
$\{\psi_{n,n+1} \colon \cD_n \to \cD_{n+1} \mid n \geq k_G\}$. 
\end{prop}
 
\subsection{Proof of Theorem~\ref{thm-main2}}\label{subsec-theorem2} 
We will prove Theorem~\ref{thm-main2} for $N = 3$, and the proof for $N \geq 3$ is verbatim the same.
Let $G \subset {\bf SL}_3(\mZ)$ be a torsion-free subgroup of finite index in the $3 \times 3$ integer matrices. Then $G$ is necessarily finitely generated. Wild actions of $G$ are constructed by the method in Section \ref{subsec-Lwild}, and we now have to show that there is an uncountable number of distinct homeomorphism types of the suspensions of such actions.

Let  $k_G \geq 2$ be chosen so that for $\cP(G) = \{p_i \mid i \geq k_G\}$  we have 
 \begin{equation}\label{eq-embedppp}
\prod_{i \geq k_G} \ {\bf SL}_3(\widehat{\mZ}_{p_i}) \subset \whG \quad \ , \quad {\rm and ~ set} ~ \whH = \prod_{p_i \in \cP(G)} {\bf SL}_3(\mZ/p_i\mZ) \ .
\end{equation}

 For each $ i \geq k_G$ we next choose a non-trivial subgroup $A_{p_i}  \subset {\bf SL}_3(\mZ/p_i\mZ)$ with trivial core.  The idea is to choose, for each $p_i$  two 
subgroups  of ${\bf SL}_3(\mZ/p_i \mZ)$ which have   orders $p_i$ and $p_i^2$, and hence cannot be isomorphic. Set  
\begin{equation}
A^1_{p_i}  =  \left\{  \left[ \begin{array}{ccc} 1 & 0 & 0 \\ 0 & 1 & a \\ 0 & 0 & 1 \end{array}\right] \mid a \in \mZ/p_i \right\} \quad , \quad A^2_{p_i} =   \left\{  \left[ \begin{array}{ccc} 1 & 0 & b \\ 0 & 1 & a \\ 0 & 0 & 1 \end{array}\right] \mid a, b \in \mZ/p_i \right\} \ .
\end{equation}
Then both $A^1_{p_i}$ and $A^2_{p_i}$ have trivial cores in ${\bf SL}_3(\mZ/p_i\mZ)$.

Let  $\ds \cI \in \Sigma_{k_G} \equiv \prod_{i \geq k_G} \ \{1,2\}$, so $\cI = (s_{k_G}, s_{k_G +1}, \ldots )$ where each $s_i \in \{1,2\}$. Now define
\begin{equation}
\cD^{\cI} = \prod_{i \geq k_G} \ A^{s_i}_{p_i} \subset \whH \ .  
\end{equation}
Then set $X^{\cI} = \whH/\cD^{\cI}$ and let  $(X^{\cI}, G, \Phi_0^{\cI})$ be the associated minimal Cantor action,   constructed as in Section \ref{subsec-Lwild}.

\begin{lemma}\label{lem-tail-equiv}
Let $\cI_1 = (s_{k_G}, s_{k_G +1},\ldots)$ and $\cI_2 =(t_{k_G}, t_{k_G+1},\ldots)$ be sequences in $\Sigma_{k_G}$. Then the sequences of surjective homomorphisms $\{\psi^1_{n,n+1}:\cD^{\cI_1}_n \to \cD^{\cI_1}_{n+1}\}$ and $\{\psi^2_{n,n+1}:\cD^{\cI_2}_n \to \cD^{\cI_2}_{n+1}\}$ are tail equivalent if and only if there exists $n \geq k_G$ such that for all $i \geq n$ we have $s_i = t_i$.
\end{lemma}

\proof By \eqref{restricted-discr}, we have 
$\ds  \cD^{\cI_1}_n \cong \prod_{i \geq n} \ A^{s_i}_{p_i} \quad , \quad \cD^{\cI_2}_n \cong \prod_{i \geq n} \ A^{t_i}_{p_i}$. 

Thus,  if $\cI_1$ and $\cI_2$ satisfy the conditions of the lemma for some $n \geq k_G$, then $\cD^{\cI_1}_i$ and $\cD^{\cI_2}_i$ are isomorphic for all $i \geq n$, and so the sequences of surjective homomorphisms of discriminant groups are tail equivalent.

Now assume that the sequences associated to $\cI_1$ and $\cI_2$ are tail equivalent.   Then there exists sequences $\{n_i\}_{ i\geq k_G}$ and $\{m_i\}_{i \geq k_G}$ and a sequence of surjective homomorphisms $\tau_{n_i, m_i}$ and $\kappa_{m_i,n_{i+1}}$,
  $$\xymatrix{ \cD^{\cI_1}_{n_1}\ar[rr]_{\tau_{n_1, m_1}}  && \cD^{\cI_2}_{m_1} \ar[rr]_{\kappa_{m_1,n_2}} &&  \cD^{\cI_1}_{n_2}  \ar[rr]_{\tau_{n_2, m_2}}  && \cD^{\cI_2}_{m_2}  \ar[rr] && \cdots }$$
 such that 
 $$\psi^1_{n_{i},n_{i+1}}= \tau_{n_i, m_i} \circ \kappa_{m_i,n_{i+1}} \quad , \quad \psi^2_{m_{i},m_{i+1}} =   \kappa_{m_i,n_{i+1}} \circ \tau_{n_{i+1},m_{i+1}} \ .$$  
 Then set  $n =  \max\{ m_i,n_{i+1}\}$. For $j \geq n$, consider a subgroup of $\cD_{m_i}^{\cI_2}$ isomorphic to $A_{p_j}^{t_j}$. Since $p_j$ is a prime, it's image under $\kappa_{m_i,n_{i+1}}$ is the subgroup of $\cD_{n_{i+1}}^{\cI_1}$, isomorphic to $A_{p_j}^{s_j}$, and it's preimage under $\tau_{n_i,m_i}$ is a subgroup of  $\cD_{n_i}^{\cI_1}$, isomorphic to $A_{p_j}^{s_j}$.
 
 If $t_j = 1$, then $s_j = 1$, as a subgroup of order $p_j$ cannot be mapped surjectively on a subgroup of order $p_j^2$. If $t_j=2$, then $s_j=2$ as the preimage of $A_{p_j}^{t_j}$ in $\cD^{\cI_1}_{n_i}$ must contain at least the same number of elements  as $A_{p_j}^{t_j}$. Thus the systems must be tail equivalent.
 \endproof

 The number of such tail equivalence classes in $\Sigma_2$ is   uncountable. Thus we conclude  that 
  there exists uncountably many distinct homeomorphism types of  weak solenoids which are wild, all with the same base manifold $M_0$ having fundamental group $G$.    This completes the proof of Theorem~\ref{thm-main2}.

We conclude with some remarks about the above construction of examples $(X^{\cI}, G, \Phi_0^{\cI})$ of minimal Cantor actions. 
First, the choice of the space $\whH$ and the closed subgroups $\cD^{\cI}$ with trivial core subgroups is somewhat obvious, given the fact that    
there exists a cofinite subset $\cP(G)$ of primes so that \eqref{eq-embedp} holds for any choice of a subgroup of finite index $G \subset   {\bf SL}_N(\mZ)$, when $N \geq 3$. This latter fact relies on the deep result, the solution of the Congruence Subgroup Problem by Bass, Milnor and Serre \cite{BMS1967}, as explained by Lubotzky in the proof of \cite[Proposition~2]{Lubotzky1993}. There are many variations of this choice families of   discriminant groups $\cD^{\cI}$, especially as $N\geq 3$ increases. The remarkable fact is that there is a finitely generated group $G$ which then acts minimally on the product space $\whH$.

Another  method to construct the Cantor space $\whH$ is to choose a countably infinite collection  of finite simple groups, $\{H_i \mid i \geq 1\}$,  which are pairwise non-isomorphic. For example, one can let $H_i = Alt(i)$ be the alternating group on $i$ symbols, for $i \geq 5$. Then for each such finite group $H_i$, we can then choose any proper, non-trivial subgroup $A_i \subset  H_i$ and it will have trivial core. 
Moreover, one consequence of the Classification Theorem for finite simple groups is that   each $H_i$  can be generated by two elements, as shown in  \cite{AG1984,Guralnick1986}. Set $\ds \whH = \prod_{i \geq 1} \ H_i$  and $\ds \cD = \prod_{i \geq 1} \ A_i$. Solutions to the following problem   then yields further  classes of examples of wild actions with uncountably many distinct asymptotic discriminant invariants.

\begin{prob}
Characterize the  countably infinite collections  of finite simple groups $\{H_i \mid i \geq 1\}$  which are pairwise non-isomorphic, such that their product $\whH$ contains a finitely generated  dense  subgroup $H$ for the Tychonoff topology on $\whH$?
\end{prob}
Then any such product group $\cD$ as above  can be realized as the discriminant of some action by a finitely generated group $H$ by the result     \cite[Theorem~10.8]{DHL2016c}, but it is unknown if these actions must be stable, or may be unstable.

Finally, we observe that in all of the known  examples of equicontinuous minimal Cantor actions $(X, G, \Phi)$ which are wild,  the action $\Phi \colon G \to Homeo(X)$ has infinitely-generated kernel.

\begin{quest}
Are there equicontinuous minimal Cantor actions $(X, G, \Phi)$ which are wild, and the kernel of $\Phi \colon G \to Homeo(X)$ is finitely-generated?
\end{quest}

\vfill
\eject


\end{document}